\documentclass[12pt,a4paper]{amsart}
\usepackage[hmargin=2.5cm,vmargin=2.5cm]{geometry}
\usepackage{url, hyperref}
\usepackage{amssymb,amsmath,amsthm}
\usepackage{paralist,enumerate,enumitem}
\usepackage{comment}

\usepackage{multirow}
\usepackage{array}

\usepackage[linecolor=red!40!white, backgroundcolor=blue!05!white,bordercolor=blue!40!white,textsize=small]{todonotes}

\tikzset{notestyleraw/.append style={align=justify}}

\newcommand{\mycomment}[1]{}

\newtheorem{theorem}{Theorem}

\newtheorem{proposition}{Proposition}
\newtheorem{remark}{Remark}

\newtheorem{hypx}{Hypothesis}
\newenvironment{hyp}{\hypx[H\textsubscript{\thehypx}]}{\endhypx}
\newcommand{\hypref}[1]{\textup{(H\textsubscript{\ref{#1}})}}

\newcommand{\md}{\mathrm{d}}
\renewcommand{\SS}{\mathbb{S}}
\newcommand{\RR}{\mathbb{R}}
\usepackage{bbm}

\newcommand{\ie}{\textit{i.e.}~}

\newcommand{\intern}{\textnormal{int}}
\newcommand{\kin}{\textnormal{kin}}
\begin{document}

\title[Compactness of linearized Boltzmann operators for polyatomic gases]{Compactness of linearized Boltzmann operators \\ for polyatomic gases
}

\author[N. Bernhoff]{Niclas Bernhoff}
\address{NB: Department of Mathematics and Computer Science, Karlstad University, Universitetsgatan 2, 65188 Karlstad, Sweden}
\email{niclas.bernhoff@kau.se}

\author[L. Boudin]{Laurent Boudin}
\address{LB: Sorbonne Universit\'e, CNRS, Universit\'e de Paris, Laboratoire Jacques-Louis Lions (LJLL), F-75005 Paris, France}
\email{laurent.boudin@sorbonne-universite.fr}

\author[M.  \v{C}oli\'{c}]{Milana  \v Coli\'{c}}
\address{M\v{C}: Department of Mathematics and Informatics, Faculty of Sciences, University of Novi Sad, Trg Dositeja Obradovi\'ca 4, 21000 Novi Sad, Serbia  } 
\email{milana.colic@dmi.uns.ac.rs}

\author[B. Grec]{B\'er\'enice Grec}
\address{BG: Université Paris Cité, CNRS, MAP5, F-75006 Paris}
\email{berenice.grec@u-paris.fr}

\thanks{This contribution is based upon work from COST Action CA18232 MAT-DYN-NET, supported by COST (European Cooperation in Science and Technology).}


\maketitle

\begin{abstract}
    In this article, we recall various existing kinetic models of non-reactive polyatomic gases. We also review the results, all recently obtained, about the compactness of the associated linearized Boltzmann operator, and briefly investigate the mixture case. Eventually, with a specific collision kernel, we present a discussion about the physical relevance of the assumptions used to obtain the compactness results.
\end{abstract}

\section{Introduction}

The compactness properties of the linearized Boltzmann collision operator $\mathcal{L}$ play a central role in the study of fluid dynamical approximations and convergence towards equilibrium for solutions of the corresponding Boltzmann equation. This linearized operator is obtained by considering a perturbation of an equilibrium, which is a function characterized in the so-called $H$-theorem as a Maxwell distribution. The operator $\mathcal{L}$ appears as the sum of a negative multiplication operator involving the collision frequency $\nu>0$, and an integral operator $K$ which is known to satisfy the expected compactness property, the proof of which we review in this paper when dealing with polyatomic gases. 

It all started with the founding paper \cite{Hilbert} by Hilbert. For more than a century now, many contributions followed this line of research, most of them are discussed in \cite{BS-16}. A major step was performed by Grad \cite{Gr-63} who obtained the compactness in the single-species monatomic case. An extension to mixtures of monatomic gases, requiring another approach to deal with disparate molecular masses, was taken care of in \cite{BGPS-13}. 

In recent years, significant efforts were made to extend the known results to physically more realistic gases,  such as the ones which involve polyatomic species. The effect of the internal structure of a polyatomic molecule is reflected on the energy conservation law during a collision. Usually, an additional internal energy variable is introduced, and this variable can be of two kinds: continuous, as in \cite{bou-des-let-per, des-mon-sal, bar-bis-bru-des}, or discrete \cite{WangChangUhlBoer, Groppi-Spiga, Gio, Kusto-book}. Recently, a generalized framework unifying these two approaches was proposed \cite{BorBisGro-22,BisBorGro-24}.

The new internal energy variable implies more difficulties to tackle the analysis of compactness properties. In the last two years, several contributions arose. They of course differ by the way the internal energy is considered (continuous or discrete), but also by the assumptions which are performed, and by the proof methodology. Those differences are the main motivation of the present work. In particular, we aim to review  those contributions by presenting them in a unified way, in terms of models and assumptions. 

\medskip

It is known that the compactness property of $K$ implies self-adjointness of the linearized operator $\mathcal{L}$, as the sum of two self-adjoint operators $K$ and $\nu \mathrm{Id}$, with $K$ bounded.
Since the set of Fredholm operators is closed under addition with compact operators, the compactness property of $K$ also implies Fredholmness of $\mathcal L$ if the collision frequency is coercive. 
Fredholmness has been considered for hard-sphere like, and hard potential like, collision kernels, in the case of polyatomic single species \cite{Be-23a, Be-23b,BST-24a} and for mixtures of monatomic and polyatomic species \cite{Be-24a, Be-24b,BST-24c}, possibly including chemical reactions \cite{Be-24c}. 

The domain of the linearized operator $\mathcal{L}$ is the same as the domain of the collision frequency, due the compactness of the remaining part of the linearized operator. In view of this, the  domain of the linearized operator $\mathcal{L}$ for hard-sphere like, and hard potential like, collision kernels is obtained by equivalent upper and lower bounds of the collision frequency. This has been considered for polyatomic single species \cite{Be-23b}, with improved upper bounds in \cite{DL-23} for continuous internal states, for discrete internal states in \cite{Be-23a}, for mixtures of monatomic and polyatomic species for discrete \cite{Be-24a} and continuous \cite{Be-24b} internal energy states, and including chemical reactions in \cite{Be-24c}.

\medskip

The paper is organized as follows. In Section~\ref{sec:models}, we bring unified ideas and notations for modelling of polyatomic gases in different settings: (i) continuous internal energy with the Borgnakke-Larsen procedure, (ii) resonant collisions,  (iii) discrete internal energy. Moreover, we extend  the presentation to gas mixtures composed of monatomic and polyatomic gases, both within the continuous and discrete approaches for the internal energy. Using this unified viewpoint, we write the Boltzmann equation and the corresponding linearized Boltzmann collision operator. Then, in Section~\ref{s:compactness}, we list the compactness results for the operator $K$ together with the associated hypotheses, obtained in various works \cite{Be-23a, Be-23b, BBS-23, borsoni-phd, Be-24a, Be-24b, BST-24a, BST-24c, shahine-phd}. In order to compare the main ideas of the proofs, we first focus on the single species case in Section~\ref{s: main ideas}, and then on the other polyatomic models in Section~\ref{s: main ideas other}. Finally, in Section~\ref{sec:discussion}, we discuss different hypotheses and their relevance with respect to experimental data, for the polyatomic single species case.

\section{Polyatomic gas modelling}\label{sec:models}

In this section, we describe various models of polyatomic gases existing in the literature. In some cases, we also explain the extension to multi-species gaseous mixtures. 

\medskip

Let us first briefly recall the description of a single monatomic gas. Consider two monatomic molecules of mass $m$ undergoing an elastic collision process. Denote their pre-collisional velocities $v$, $v_*$, which become $v'$, $v'_*$ after collision. The microscopic momentum and energy are conserved through the process, expressed as follows
\begin{align}
m v + m v_* &= m v'+ mv'_*,  \phantom{\frac 12}\label{e:monomicromomentcons}\\
\frac{1}2 m|v|^2 + \frac{1}2 m |v_*|^2 &= \frac{1}2 m|v'|^2 + \frac{1}2 m|v'_*|^2.   \label{e:monomicroenercons}
\end{align}
Introducing a parameter $\sigma\in\SS^2$, one obtains the collision rules, giving the expression of the post-collisional velocities
\begin{equation} \label{e:paramvelocmono}
v'=\frac12 (v+v_*) + \frac12 |v-v_*|\sigma, \qquad v_*'=\frac12 (v+v_*) - \frac12 |v-v_*|\sigma,
\end{equation}
in terms of pre-collisional ones. For (much) more details on the monatomic case, the reader may refer, for instance, to \cite{villani}.

\medskip

At the macroscopic level, the main distinction between monatomic and polyatomic gases can be observed when the specific internal energy of the gas is considered. Let us denote by $\hat{e}$ its dimensionless form. Namely, for a monatomic gas, it is known that the dimensionless specific heat at constant volume is  $\hat{c}_v=3/2$ and consequently $\frac{\md\hat{e}}{\md T}= 3/2$. However, for  thermally perfect (non-polytropic) gases, the specific heat $\hat{c}_v$ is a temperature-dependent quantity and, in general, it is only known through its definition $\hat{c}_v(T)=\frac{\md\hat{e}}{\md T}$. Even if $\hat{c}_v(T)$ is assumed to be constant with respect to the temperature, which corresponds to calorically perfect or polytropic gases, it is experimentally observed that, at room temperature, $\hat{c}_v > 3/2$ \cite{Aoki, MPC-SS-non-poly}. Hence, such a behavior of polyatomic gases is a consequence of more complex collisions than the ones for monatomic gases, and in particular, \eqref{e:monomicroenercons} has to be reconsidered. 

\medskip

Two approaches to model polyatomic gases have been developed in parallel. They share the same idea: associate an internal energy to a polyatomic molecule and consequently rewrite \eqref{e:monomicroenercons}. The main difference lies in the form of the internal energy they use, either a \emph{continuous} or \emph{discrete} one.

The continuous approach introduces a continuous internal energy $I\in \RR_+$ as an additional argument of the distribution function, and a nonnegative function $\varphi$ of $I$ which becomes a parameter of the model to capture a proper form of the specific internal energy $\hat{e}$. First, \cite{bou-des-let-per} proposed a power-law form of $\varphi$, that is, for any $I\geq0$,
\begin{equation}\label{power law phi}
    \varphi(I)=I^{\delta/2-1},
\end{equation}
where  $\delta>0$ is related to the number of internal degrees of freedom of the molecules. More specifically, for polytropic gases, $\delta$ is  constant and related to measurements \cite{MPC-Dj-T-O} via 
\begin{equation}\label{delta cv}
    \delta=2 \hat{c}_v -3.
\end{equation}
 The model is accurate for diatomic gases (for instance $\mbox{N}_2$,  $\mbox{O}_2$, CO, $\mbox{H}_2$) for temperatures close to the room temperature. Later, \cite{des-mon-sal} presented the model without prescribing a specific form of $\varphi$. This general form is supposed to allow a more general macroscopic internal-energy law and eventually capture non-polytropic gases. 

With the discrete approach \cite{WangChangUhlBoer,Gio,Groppi-Spiga}, the internal energy can only take a finite number of given values $I^{(1)}$, ..., $I^{(N_\intern)}$, with $N_\intern \ge 2$. Then the parameters $\varphi^{(1)}$, ..., $\varphi^{(N_\intern)}$  account for the degeneracy of the different energy levels \cite{Gio}. 
The degeneracy $\varphi^{(k)}\in \RR_+$ of the internal energy $I^{(k)}$ corresponds to the number of different states that give rise to the same specific internal energy $I^{(k)}$.

\subsection{Borgnakke-Larsen procedure} \label{ss:BL}

Originally presented in \cite{bor-lar-75}, this model was the first one studied in the literature when dealing with a continuous internal energy variable in a kinetic setting. It can be found, written in various forms, for instance in \cite{bou-des-let-per, des-mon-sal}. 

Let us consider two colliding  polyatomic molecules of mass $m$ with respective velocities and internal energies $(v, I)$ and $(v_*, I_*)$, changing to  $(v', I')$ and $(v'_*, I'_*)$. The microscopic momentum and total energy are conserved as follows
\begin{align} 
m v + m v_* &= m v'+ mv'_*, \phantom{\frac 12}\label{microCL_v}\\ 
\frac{m}{2} |v|^2 + I + \frac{m}{2} |v_*|^2 + I_* &= \frac{m}{2} |v'|^2 + I'+ \frac{m}{2} |v'_*|^2 + I'_*.\label{microCL_v^2}
\end{align}
In the center-of-mass reference frame, using the relative velocities $V=v-v_*$, $V'=v'-v'_*$, the energy conservation can equivalently be rewritten as
\begin{equation}  \label{energy}
E:=\frac{m}{4} \left|V\right|^2 + I + I_* = \frac{m}{4} \left|V'\right|^2 + I'+ I'_*,
\end{equation}
which defines the total energy $E$ in the center-of-mass frame.
The Borgnakke-Larsen procedure splits $E$ into  kinetic and internal energy contributions using a parameter $R\in[0,1]$, as
\begin{equation}  \label{totalenergypartition}
	RE=\frac{m}{4} \left|V'\right|^2, \qquad (1-R) E = I' + I_*',
\end{equation}
and then  associates those energies to the post-collisional velocities and internal energies with the help of parameters $\sigma \in \mathbb{S}^2$ and $r\in[0,1]$, with
\begin{gather} 
v'= \frac{v+v_*}{2} + \sqrt{\frac{R E}{m}} \sigma, \qquad v'_* = \frac{v+v_*}{2} - \sqrt{\frac{R E}{m}}\sigma, \label{colllaws_v}\\
I'= r (1-R) E, \qquad I'_* = (1-r) (1-R) E.  \label{colllaws_I}
\end{gather}
The collision kernel  $B\geq 0$  is assumed to satisfy symmetry properties reflecting an interchange of colliding molecules, as well as microreversibility assumptions corresponding to the pre/post-collision change, which means that
\begin{equation}\label{assumptions_B_BL}
\begin{split}
B(v,v_*,I,I_*,r,R,\sigma)&=B(v_*,v,I_*,I,r,R,\sigma),\\
B(v,v_*,I,I_*,r,R,\sigma)&=B\left(v',v'_*,I',I'_*,r',R',\sigma'\right),
\end{split}
\end{equation}
with 
\begin{equation} \label{e:defRrprime}
R'=\frac{m |V|^2}{4 E}, \qquad r'= \frac{I}{I+I_*} = \frac{I}{(1-R')E}, \qquad \sigma'=\frac{V}{|V|}.
\end{equation}
The associated Boltzmann collision operator can then be defined, for any measurable function $f$ for which it makes sense, and for almost every $(v,I) \in \RR^3\times\RR_+$, by
\begin{multline}  \label{collisionoperatorBLgeneral}
Q(f,g)(v,I) = \int_{\RR^3\times \RR_+\times (0,1)^2\times \SS^2 } \left(f'\, g'_*  \frac{\varphi(I) \varphi(I_*) }{
\varphi(I') \varphi(I'_*)}   - f \, g_* \right) \\\phantom{\int}
\times \frac{\tilde{B}(v,v_*,I,I_*,r,R,\sigma)}{\varphi(I) \varphi(I_*)} \,  (1-R) \sqrt{R}\, \md v_* \, \md I_* \,\md R \, \md r \, \md \sigma,
\end{multline}
with $\tilde{B}$ satisfying \eqref{assumptions_B_BL}, and where we used the standard notations $f'=f(v',I')$, $g'_*=g(v'_*,I'_*)$, $f=f(v,I)$, $g_*=g(v_*,I_*)$, the prime quantities $v'$, $v_*'$, $I'$ and $I_*'$ being defined by \eqref{colllaws_v}--\eqref{colllaws_I}. Results that are of interest for this review use the specific form \eqref{power law phi} of the weight factor, namely $\varphi(I)=I^{\delta/2-1}$, $\delta>0$. In this case, one can use another collision kernel defined as
\begin{equation*}
    B(v, v_*, I, I_*, r, R, \sigma) = \frac{\tilde{B}(v,v_*,I,I_*,r,R,\sigma)}{(I I_*)^{\delta/2-1} (r(1-r))^{\delta/2-1} (1-R)^{{\delta-2}}}
    =\frac{\tilde{B}(v,v_*,I,I_*,r,R,\sigma)E^{\delta-2}}{\varphi(I) \varphi(I_*) \varphi(I')\varphi(I_*')},
\end{equation*}
which clearly still satisfies assumption \eqref{assumptions_B_BL}. Then the collision operator reads,  for almost every $(v,I) \in \RR^3 \times \RR_+$,
\begin{multline}  \label{collisionoperator}
Q(f,g)(v,I) = \int_{\RR^3\times \RR_+\times (0,1)^2\times \SS^2 } \left(f'\, g'_* \left(\frac{I \, I_* }{%
I'\, I'_*} \right)^{\delta/2-1} - f \, g_* \right) \\
\phantom{\int}\times B(v,v_*,I,I_*,r,R,\sigma) \, \left( r (1-r) \right)^{\delta/2 - 1} (1-R)^{\delta - 1} \sqrt{R}\,
\md v_* \, \md I_* \,\md R \, \md r \, \md \sigma.
\end{multline}

For convenience, we state a part of the $H$-theorem that defines the equilibrium state. 
\begin{proposition} \label{prop:htheorem_BL}
The three following properties are equivalent:
\begin{enumerate}[label=(\roman*)]
\item $\displaystyle Q(M,M) = 0$,$\displaystyle \phantom{\int}$ 
\item $\displaystyle\int_{\RR^3 \times \RR_+} Q(M,M) (v,I)  \, \log (M (v,I) I^{1-\delta/2}) \, \md v \, \md I= 0$, 
\item there exist $n \geq 0$, $u \in \RR^3$ and  $T > 0$ such that, for almost every $(v,I) \in \RR^3 \times \RR_+$, 
\begin{equation} \label{e:defMaxwBL}
 M(v,I) = \frac{n}{(k_B T)^{\delta/2} \Gamma(\delta/2)} \left( \frac{m}{2 \pi k_B T}\right)^{3/2} I^{\delta/2-1} \exp \left(- \frac{m|v-u|^2}{2 k_B T}  - \frac{I}{k_B T} \right),
\end{equation}
where $\Gamma$ represents the usual Gamma function. 
\end{enumerate}
\end{proposition}

\subsection{Polyatomic gas with resonant collisions}

A resonant behaviour can be observed, for instance, in the collisions between selectively excited $\textnormal{CO}_2$ molecules \cite{ree-et-al-06}. In that case, the microscopic internal and kinetic energies are separately conserved during the collisional process. 

Consider two resonant-colliding polyatomic molecules of mass $m$ with velocities and internal energies  $(v, I)$ and $(v_*, I_*)$, which change into $(v', I')$ and $(v'_*, I'_*)$ due to the collision process. In the resonant kinetic model introduced in \cite{bou-ros-sal-22}, we have the following microscopic momentum and energy conservations \begin{align}
m v + m v_* &= m v'+ mv'_*,  \phantom{\int}\label{e:resmicromomentcons}\\
\frac{m}2 |v|^2 + \frac{m}2  |v_*|^2 &= \frac{m}2 |v'|^2 + \frac{m}2 |v'_*|^2, \label{e:resmicroenercons_v}\\
I + I_* &=I'+ I'_*. \phantom{\int} \label{e:resmicroenercons_I}
\end{align}
Due to the conservations \eqref{e:resmicromomentcons}--\eqref{e:resmicroenercons_v}, which are separate, we inherit the natural monatomic parametrization \eqref{e:paramvelocmono} of the velocities. Then, in \cite{bou-ros-sal-22}, the internal energies are parametrized using the same idea as in the Borgnakke-Larsen model recalled in the previous subsection. More precisely, the authors introduce a parameter allowing to distribute the conserved internal energy between the prime internal energies of both molecules. 
\begin{equation} \label{e:defI*'reson}
I_*'=I+I_*-I'.
\end{equation}
Let us describe the model with this last parametrization. We consider a collision kernel $B\geq 0$ which is required to be symmetric, \ie
\begin{equation} \label{e:Bsymm}
B(v,v_*,I,I_*,I',\sigma) = B(v_*,v,I_*,I,I'_*,\sigma), 
\end{equation}
and to satisfy a microreversibility property, that is  
\begin{equation} \label{e:Bmicrorev}
B(v,v_*,I,I_*,I',\sigma) = B\left(v',v'_*,I',I'_*,I,\sigma'\right),
\end{equation}
where $V=v-v_*$ denotes again the relative velocity. The associated invariant measure is
$$B(v,v_*,I,I_*,I',\sigma)\,\mathbf{1}_{\left[0, I+I_* \right]}(I')\,\varphi(I)\,\varphi(I_*)\,\varphi(I')\,\varphi(I+I_*-I')\,\md I'\,\md I_*\,\md I\,\md \sigma\,\md v_*\,\md v,$$
where $\mathbf{1}_{\left[0, I+I_* \right]}$ is the characteristic function of $[0,I+I_*]$. 
Note that \eqref{e:Bmicrorev} is the microreversibility condition written in both \cite{BBS-23} and \cite{borsoni-phd}, but the measure was inaccurately written in \cite{BBS-23} and corrected in \cite{borsoni-phd, resonant-corr}. 

\begin{remark}
A resonant collision kernel is related to a collision kernel similar to the ones used in Subsection~\ref{ss:BL} in the following way. Consider a collision kernel function $
\hat B(v,v_*,I,I_*,r,R,\sigma)$ satisfying the symmetry/microreversiblity conditions \eqref{assumptions_B_BL}.
Following \cite{borsoni-phd}, and keeping the notation $E$ for the center-of-mass energy, obviously conserved during the collision process, one can build a resonant collision kernel $B$ from $\hat B$ through
\begin{equation}\label{equiv-B-Bhat}
B(v,v_*,I,I_*,I',\sigma)=\hat B\left(v,v_*,I,I_*,\frac{I'}{I+I_*},1-\frac{I+I_*}E,\sigma\right) \, E \,\mathbf{1}_{\left[0, I+I_* \right]}(I').
\end{equation}
Note that the characteristic function can be dropped (it is only added here for the sake of clarity) since we can assume that $\hat B$ is zero when computed for values of parameters outside $[0,1]$. \end{remark}

This allows to define the Boltzmann collision operator $Q$ associated to the resonant model. For any measurable function $f$ for which it makes sense, we can write, for almost every $(v,I) \in \RR^3 \times \RR_+$,
\begin{multline} \label{e:defQreson}
    Q(f,g)(v,I) = \int_{\RR^3 \times (\RR_+)^2 \times \SS^2} \left( f' g'_*\,\frac{\varphi(I)\varphi(I_*)}{\varphi(I')\varphi (I+I_*-I')} - f g_* \right) \\
    \phantom{\int}\times B(v,v_*,I,I_*,I',\sigma) \,\mathbf{1}_{\left[0, I+I_* \right]}(I') \,\frac{\varphi(I') \, \varphi(I+I_*-I')}{\Psi_{\text{res}}(I+I_*)}  \, \md v_* \, \md I_* \, \md I'  \, \md \sigma,
\end{multline}
where we used again the standard notations $f'=f(v',I')$, $g'_*=g(v'_*,I'_*)$, $f=f(v,I)$, $g_*=g(v_*,I_*)$, with the prime quantities $v'$, $v_*'$ and $I_*'$ defined by \eqref{e:paramvelocmono} and \eqref{e:defI*'reson}, and, for any $Z\ge 0$, 
$$\Psi_{\text{res}}(Z)=\int_0^Z \varphi (I')\,\varphi (Z-I')\,\md I'.$$ 

The equilibria are characterized thanks to the following part of the $H$-theorem proved in \cite{BBS-23}.
\begin{proposition} \label{prop:htheorem_reson}
The three following properties are equivalent:
\begin{enumerate}[label=(\roman*)]
\item $\displaystyle Q(M,M) = 0$,$\displaystyle \phantom{\int}$ 
\item $\displaystyle\int_{\RR^3 \times \RR_+} Q(M,M) (v,I) \, \log \left[ \frac{M (v,I)}{ \varphi(I)} \right] \, \md v \, \md I= 0$, 
\item there exist $n \geq 0$, $u \in \RR^3$ and $T_\kin$, $T_\intern > 0$ such that, for almost every $v$ and $I$, 
\begin{equation} \label{e:defMaxwreson}
 M(v,I) = \frac{n}{q(T_\intern)} \left( \frac{m}{2 \pi k_B T_\kin}\right)^{3/2} \varphi(I) \exp \left(- \frac{m|v-u|^2}{2 k_B T_\kin}  - \frac{I}{k_B T_\intern} \right),
\end{equation}
\end{enumerate}
where the internal energy partition function $q$ is defined, for any $T>0$, by 
\begin{equation} \label{q_cont}
q(T) = \int_{\RR_+} \exp\left( -\frac I{k_B T}\right) \varphi(I)\,\md I.
\end{equation}
\end{proposition}
Note that \eqref{e:defMaxwreson} is a product of two Gibbs distributions, which can also be named a Maxwell distribution similarly to the monatomic case, with two different temperatures $T_\kin$ and~$T_\intern$. 

\begin{remark}
In \cite{BBS-23}, the authors provide a framework for $\varphi$ of the following form: there exist $\beta_1$, $\beta_2\ge0$, and $C$, $C'>0$ such that 
$$C I^{\beta_1}\le\varphi(I)\le C' I^{\beta_1}, \qquad I\in (0,1),$$
and that, for any $a>0$, there exists $C_a>0$ such that
$$C_a I^{\beta_2-a}\le\varphi(I)\le C' I^{\beta_2}, \qquad I\ge1.$$
\end{remark}
For the sake of upcoming discussions, we rewrite the collision operator \eqref{e:defQreson} and equilibrium distributions \eqref{e:defMaxwreson} by explicitly using \eqref{power law phi} for $\varphi$. First, one can compute $\Psi_{\text{res}}(Z) = Z^{\delta-1} \Gamma(\delta/2)^2/\Gamma(\delta)$, and thus \eqref{e:defQreson} reads, for almost every $v$ and $I$,
\begin{multline} \label{e:defQreson 2}
    Q(f,g)(v,I) = \frac{ \Gamma(\delta)}{\Gamma(\delta/2)^2}
    \int_{\RR^3 \times (\RR_+)^2 \times \SS^2} \left( f' g'_*\,\left(\frac{I I_*}{I'(I+I_*-I')}\right)^{\delta/2-1} - f g_* \right) \\
    \phantom{\int}\times B(v,v_*,I,I_*,I',\sigma) \,\mathbf{1}_{\left[0, I+I_* \right]}(I')
    \frac{\left[I' (I+I_*-I')\right]^{\delta/2-1}}  {(I+I_*)^{\delta-1}} \, \md v_* \, \md I_* \, \md I'  \, \md \sigma.
\end{multline}

The equilibrium \eqref{e:defMaxwreson} becomes, for almost every $v$ and $I$,
\begin{equation} \label{e:defMaxwresonpowerlaw}
 M(v,I) = \frac{n}{\Gamma(\delta/2)} \left( \frac{m}{2\pi k_BT_\kin}\right)^{3/2}(k_B T_\intern)^{-\delta/2}\, I^{\delta/2-1} \exp \left(- \frac{m|v-u|^2}{2 k_B T_\kin}  - \frac{I}{k_B T_\intern} \right).
\end{equation}

\subsection{Polyatomic gas with discrete internal energies}\label{sec:discrete-single}

Let us now describe another approach to model the degrees of freedom associated to internal energy for a single polyatomic gas. Instead of considering a continuous internal energy variable, one can introduce $N_\intern\ge 2$ different internal energies $I^{(1)}$, ..., $I^{(N_\intern)}\in \mathbb{R}_{+}$, see for instance \cite{WangChangUhlBoer,Gio,Groppi-Spiga}.
Each collision can then be represented by two pre-collisional pairs and two corresponding post-collisional pairs, respectively indexed by $(k,\ell)$ and $(k',\ell')$, with $k$, $\ell$, $k'$, $\ell' \in \left\{1,...,N_\intern\right\}$. More precisely, consider two colliding molecules of mass $m$ with  velocities and internal energies $( v,I^{(k)}) $ and $( v_* ,I^{(\ell)})$, changing into $( v',I^{(k')}) $ and $( v'_*,I^{(\ell')}) $, and define the internal energy gap
\begin{equation*}
\Delta I^{(k\ell,k'\ell')} =I^{(k')}+I^{(\ell')}-I^{(k)}-I^{(\ell)}.
\end{equation*} 
The microscopic momentum and energy conservations are now written as
\begin{align}
m v + m v_* &= m v'+ mv'_*, \\    
\frac{m}{2}|v|^2+\frac{m}{2}|v_*|^2+I^{(k)}+I^{(\ell)} &=\frac{m}{2}|v'|^2+\frac{m}{2}|v'_*|^2+I^{(k')}+I^{(\ell')}.
\end{align}
Equivalently, the conservations can be rewritten, defining the total energy $E^{(k\ell)}$ in the center-of-mass reference frame, as
\begin{equation*}
E^{(k\ell)}:=\dfrac{m}{4}\left|V\right|^2+I^{(k)}+I^{(\ell)}=\dfrac{m}{4}\left|V'\right|^2+I^{(k')}+I^{(\ell')}.
\end{equation*}%
The post-collisional velocities are then given by
\begin{equation*}
v'=\dfrac{v+v_*}{2}+\sqrt{|V|^2-\dfrac{4\Delta I^{(k\ell,k'\ell')} }{m}}\dfrac{\sigma}{2}, \qquad v'_*=\dfrac{v+v_*}{2}-\sqrt{|V|^2-\dfrac{4\Delta I^{(k\ell,k'\ell')}}{m}}\dfrac{\sigma}{2}.
\end{equation*}%
The nonnegative collision kernel $B$ is again assumed to satisfy symmetry and microreversibility relations as follows
\begin{gather*}
B( v,v_*,I^{(k)},I^{(\ell)},I^{(k')},I^{(\ell')},\sigma) 
=B( v_*,v,I^{(\ell)},I^{(k)},I^{(k')},I^{(\ell')},\sigma )\\
=B( v,v_*,I^{(k)},I^{(\ell)},I^{(\ell')},I^{(k')},-\sigma ) , \\
 B( v,v_*,I^{(k)},I^{(\ell)},I^{(k')},I^{(\ell')},\sigma) 
 = B\left( v',v'_*,I^{(k')},I^{(\ell')},I^{(k)},I^{(\ell)} ,\sigma'\right).
\end{gather*}

In this discrete internal energy case, the distribution function is studied under the form $f=( f^{(1)},\dots,f^{(N_\intern)}) $, where each component $f^{(k)} =f^{(k)}(t,x,v)$, $1\le k\le N_\intern$, is the distribution function for particles with internal energy $I^{(k)}$. 
This allows to define the Boltzmann collision operator associated to the $k$th component of the distribution function, for almost every $v \in \RR^3$,
\begin{multline}
Q^{(k)}(f,g)(v)=\sum\limits_{\ell,k',\ell'=1}^{N_\intern}\int_{\RR^3\times \SS^2} \left( f^{(k')}(v')g^{(\ell')}(v'_*)\frac{\varphi^{(k)}\varphi^{(\ell)}}{\varphi^{(k')}\varphi^{(\ell')}}- f^{(k)} (v) g^{(\ell)}(v_*)\right)\\
\times B( v,v_*,I^{(k)},I^{(\ell)},I^{(k')},I^{(\ell')},\sigma )\varphi^{(k')}\varphi^{(\ell')}\frac{|V'|}{\left(E^{(k\ell)} \right)^{1/2}}\,\md v_*\,\md\sigma.
\label{collisionoperator_discrete}
\end{multline}%

The equilibrium states are characterized by the following part of the $H$-theorem.
\begin{proposition}
\label{prop:htheorem_DI} The three following properties are equivalent:

\begin{enumerate}[label=(\roman*)]
\item $\displaystyle Q^{(k)}(M,M)=0$ for any $k\in \left\{ 1,...,N_\intern\right\}$,

\item $\displaystyle\sum\limits_{k=1}^{N_\intern}\int_{\RR^{3}}Q^{(k)}(M,M)(v)\,%
\log \left[ \frac{M^{(k)}(v)}{ \varphi^{(k)}} \right]\,\md v\,=0$,

\item there exist $n\geq 0$, $u\in \RR^{3}$ and $T>0$ such that, for every $%
k\in \left\{ 1,...,N_\intern\right\} $ and almost every $v$,%
\begin{equation}
M^{(k)}(v)=\frac{n}{q}\left( \frac{m}{2\pi k_{B}T}\right) ^{3/2}\varphi
^{(k)}\exp \left( -\frac{m|v-u|^{2}}{2k_{B}T}-\frac{I^{(k)}}{k_{B}T}\right) ,
\label{e:defMaxwDI}
\end{equation}%
with $q=\sum\limits_{k=1}^{N_\intern}\varphi ^{(k)}\exp \left( -\dfrac{I^{(k)}}{%
k_{B}T}\right) $ being the discrete version of integral \eqref{q_cont}.
\end{enumerate}
\end{proposition}

\subsection{Polyatomic models for mixtures}
\subsubsection{Borgnakke-Larsen model for a mixture of monatomic and/or polyatomic gases} \strut\\
Consider a mixture of $N\ge 2$ monatomic or polyatomic species. Let us denote by $\mathcal{M}$ the set of indices corresponding to monatomic gases, and $\mathcal{P}$ the one corresponding to polyatomic gases, so that $\mathcal{M}\cup\mathcal{P}=\{1,\dots,N\}$. 
For any polyatomic species $i\in \mathcal P$, the parameter used in \eqref{power law phi} is denoted by $\delta_i>0$, which is related to the number of internal degrees of freedom of the molecules. 
For two colliding molecules of species $i$ and $j$, with respective masses $m_i$ and $m_j$, velocities $v$ and $v_*$, and (if polyatomic) internal energies $I$ and $I_*$, the microscopic momentum and energy conservations are now written as
\begin{eqnarray*}
m_i v + m_j v_* &=&m_i v'+m_j v_*', \\
\frac{m_i}{2}|v|^2 + \frac{m_j}{2} |v_*|^2 + I \mathbf{1}_{i\in \mathcal{P}}+I_* \mathbf{1}_{j\in \mathcal{P}}
&=& \frac{m_i}{2}|v'|^2 + \frac{m_j}{2} |v'_*|^2 +I' \mathbf{1}_{i\in \mathcal{P}} + I'_* \mathbf{1}_{j\in \mathcal{P}},
\end{eqnarray*}
or equivalently with the energy conservation law   in the center-of-mass frame
\begin{equation*}
E_{ij} = \frac{\mu_{ij}}{2} |V|^2 + I \mathbf{1}_{i \in \mathcal{P}} + I_* \mathbf{1}_{j \in \mathcal{P}}  = \frac{\mu_{ij}}{2} |V'|^2 + I' \mathbf{1}_{i \in \mathcal{P}} + I'_* \mathbf{1}_{j \in \mathcal{P}}, 
\end{equation*}
where $\mu_{ij} = \frac{m_i m_j}{m_i+m_j}$ is the reduced mass. 

Let us now detail, for each case, the explicit expressions of the collision rules with a Borgnakke-Larsen-like parametrization for polyatomic gases \cite{des-mon-sal}, as well as the ones of the corresponding collision operators \cite{bar-bis-bru-des,Alonso-Colic-Gamba}.

\paragraph{\textit{\textbf{Collision between two monatomic molecules}}}

In this case, the collision rules can be written with the same parametrization $\sigma\in\SS^2$ as in the single species case, \ie
\begin{equation*} 
v'= \frac{m_i v + m_j v_*}{m_i+m_j} + \frac{m_j}{m_i + m_j} |V| \sigma, 
\qquad v'_* = \frac{m_i v + m_j v_*}{m_i+m_j} - \frac{m_i}{m_i + m_j} |V| \sigma.
\end{equation*}
The collision kernels $B_{ij}(v, v_*, \sigma) \geq 0$ are again assumed to satisfy a symmetry and a microreversibility property as follows
\begin{equation}  \label{m-m cross}
B_{ij}(v, v_*, \sigma) = B_{ji}(v_*, v, -\sigma) = B_{ij}\left(v', v'_*,\sigma'\right).
\end{equation}
\begin{remark} \label{r:intra-inter-mono}
Due to the indistinguishability of particles of the same species, in a binary intra-species collision, we do not distinguish which of the two particles has which velocity after the collision (both outcomes are assumed to be equally probable). Therefore, either pre- or post-collisional velocities can be separately interchanged without affecting the collision kernel. Noting that $\sigma=\frac{v'-v'_*}{|v'-v'_*|}$, one obtains the symmetry properties
\begin{equation*}
B_{ii}(v, v_*, \sigma) = B_{ii}(v_*, v, \sigma) = B_{ii}(v_*, v, -\sigma) = B_{ii}(v, v_*, -\sigma).
\end{equation*}
Note that the last two equalities are also consequences of the first symmetry and the microreversibility property.

In an inter-species collision, this is not the case anymore. The two particles are now distinguishable and one has to tell apart which particle has which velocity. The interchange of the particle roles, highlighted by the interchange of indices in the collision kernels, implies the interchange of both pre- and post-collisional velocities simultaneously. Therefore, the only symmetry property which remains true is \eqref{m-m cross}.
\end{remark}

The Boltzmann collision operator is then written, for almost every $v\in\RR^3$, as
\begin{equation}  \label{m-m Q}
Q_{ij}(f,g)(v) = \int_{\RR^3\times \SS^2} \left( f'g'_* - f g_* \right) B_{ij}(v,v_*,\sigma) \, \md v_*\, \md \sigma .
\end{equation}

\paragraph{\textit{\textbf{Collision between two polyatomic molecules}}}
When at least one polyatomic species is involved, the velocity collision rules become
\begin{align}
v'&= \frac{m_i v + m_j v_*}{m_i + m_j} + \frac{m_j}{m_i + m_j} \sqrt{\frac{2 \, R \, E_{ij}}{\mu_{ij}}} \sigma, \label{e:collrules_mix_poly_v'} \\
v'_{*} &= \frac{m_i v + m_j v_*}{m_i + m_j} - \frac{m_i}{m_i + m_j} 
\sqrt{\frac{2 \, R \, E_{ij}}{\mu_{ij}}} \sigma.\label{e:collrules_mix_poly_v'*}
\end{align}
 With the usual Borgnakke-Larsen parametrization, the allocation between the two microscopic internal energies provides
\begin{equation}\label{e:collrules_mix_poly_I'}
I'=r (1-R) E_{ij}, \qquad I'_* = (1-r)(1-R) E_{ij}.
\end{equation}
Again, symmetry and microreversibility are required for the collision kernels, \ie 
\begin{equation*}
    B_{ij}(v,v_*,I,I_*,r,R,\sigma)=B_{ji}(v_*,v,I_*,I,1-r,R,-\sigma)=B_{ij}\left(v',v'_*,I',I'_*,r',R',\sigma'\right).
\end{equation*}

\begin{remark} \label{r:intra-inter-poly}
When $i\neq j$, the interchange of the particle roles implies the interchange of both the pre- and post-collisional internal energies, additionally to the interchange of both the pre- and post-collisional velocities as in Remark~\ref{r:intra-inter-mono}. Thus the symmetry property has to involve both changes $\sigma$ into $- \sigma$ and $r$ into $1-r$.
\end{remark}

The collision operator for two polyatomic molecules of species $i$, $j \in \mathcal P$ reads, for almost every $(v,I)\in\RR^3\times \RR_+$,
\begin{multline}  \label{p-p Q}
Q_{ij}(f,g)(v,I) = \int_{\RR^3\times \RR_+\times(0,1)^2\times \SS^2} \left( f' g'_*\left(\frac{I }{%
I'}\right)^{\delta_i/2-1 } \left(\frac{ I_*}{ I'_*}\right)^{\delta_j/2-1} - f g_* \right) \\
 \phantom{\int}\times B_{ij}(v, v_*, I, I_*, r, R, \sigma) \, r^{\delta_i/2-1}(1-r)^{\delta_j/2-1} \, (1-R)^{\delta_i/2 + \delta_j/2-1} \, 
\sqrt{R} \,\md v_* \, \md I_* \,\md R \, \md r \, \md \sigma.
\end{multline}

\paragraph{\textit{\textbf{Collision between one polyatomic molecule and one monatomic molecule}}}
In the case of a collision between one polyatomic and one monatomic molecule, the collision rules \eqref{e:collrules_mix_poly_v'}--\eqref{e:collrules_mix_poly_v'*} still hold. Nevertheless, other expressions are changed, since there is no parameter $r$ anymore associating the internal energy to each species.

If a polyatomic molecule of species $i\in \mathcal{P}$ collides with a monatomic molecule of species $j \in \mathcal{M}$, we have for the internal energy
\begin{equation*}
I'= (1-R) E_{ij},
\end{equation*}
the symmetry and microreversibility assumptions on $B_{ij}$ are given by
\begin{equation}
    B_{ij}(v,v_*,I,R,\sigma)=B_{ji}(v_*,v,I,R,-\sigma)=B_{ij}\left(v',v'_*,I',R',\sigma'\right),
\end{equation}
and the collision operator is written, for almost every $(v,I)\in\RR^3\times \RR_+$, as
\begin{multline}  \label{p-m Q}
Q_{ij}(f,g)(v,I) = \int_{\RR^3\times(0,1)\times \SS^2} \left( f(v',I')g(v'_*)\left(\frac{I}{I'}\right)^{\delta_i/2-1} - f(v,I)g(v_*) \right) \\
\phantom{\int}\times {B}_{ij}(v, v_*, I, R,\sigma) \, (1-R)^{\delta_i/2-1} \, \sqrt{R} \,\md v_* \,\md R  \, \md \sigma.
\end{multline}

In the opposite case, when $i\in \mathcal{M}$ and $j\in \mathcal{P}$, the post-collisional internal energy is given by
\begin{equation*}
I'_{*} = (1-R)\, E_{ij},
\end{equation*}
the assumptions on $B_{ij}$ are given by
\begin{equation}
    B_{ij}(v,v_*,I_*,R,\sigma)=B_{ji}(v_*,v,I_*,R,-\sigma)=B_{ij}\left(v',v'_*,I'_*,R',\sigma'\right),
\end{equation}
and the collision operator becomes, for almost every $v\in\RR^3$,
\begin{multline}  \label{m-p coll operator}
Q_{ij}(f,g)(v) = \int_{\RR^3\times \RR_+\times(0,1)\times \SS^2} \left( f(v')g(v'_*,I'_*)\left(\frac{I_*}{I'_*}\right)^{\delta_j/2-1} - f(v)g(v_*,I_*)
\right) \\
\phantom{\int}\times {B}_{ij}(v,v_*,I_*,R,\sigma) \, (1-R)^{\delta_j/2-1} \, \sqrt{R} \,\md v_* \, \md I_* \,\md R \, \md \sigma.
\end{multline}

In the previous equalities  \eqref{m-m Q}, \eqref{p-p Q}, \eqref{p-m Q} and \eqref{m-p coll operator}, both $f$ and $g$ are scalar functions related to the species $i$ and $j$. In order to define the $i$th component $Q_i$ of the vector collision operator, we must emphasize that its argument has to be all the scalar functions related to any species at the same time. In other words, we use the vector forms $f=(f_1,\dots,f_N)$, $g=(g_1,\dots,g_N)$, and then write, for any $i\in\{1,\dots,N\}$,
\begin{equation} \label{e:ithcollop}
Q_i(f,g) =\sum_{j=1}^N Q_{ij} (f_i,g_j).
\end{equation}

Then, as usual, the equilibria are described thanks to the following part of the $H$-theorem.
\begin{proposition} \label{prop:htheorem_MBL} 
The three following properties are equivalent:
\begin{enumerate}
\item $\displaystyle Q_i(M,M)=0$ for any $i$ such that $1\leq i \leq N$,
\item the following equality holds
\begin{multline*}
\sum\limits_{i\in \mathcal{M}}\int_{\RR^{3}}Q_{i}(M,M)(v)\,
\log M_{i}(v)\,\md v\\
+\sum\limits_{i\in \mathcal{P}}\int_{\RR^{3}\times \RR_+}Q_{i}(M,M)(v,I)\,%
\log \left ( M_{i}(v,I)I^{1-\delta_i/2}  \right )\,\md v\, \md I=0,
\end{multline*}

\item there exist $n=(n_1,\dots,n_N)\in \RR_+^N$, $u\in \RR^{3}$ and $T>0$
such that, if $i\in \mathcal{M}$,  for almost every $v$,
\begin{equation} \label{e:defMaxwMBLmono}
M_{i}(v)=n_{i}\left( \frac{m_{i}}{2\pi k_{B}T}\right)
^{3/2}\exp \left( -\frac{m_{i}|v-u|^{2}}{2k_{B}T}\right), 
\end{equation}
and if $i\in \mathcal{P}$,  for almost every $v$ and $I$,
\begin{equation} \label{e:defMaxwMBLpoly}
M_{i}(v,I)=\frac{n_{i}}{\Gamma(\delta_i/2)}\left( \frac{m_{i}}{2\pi}\right)^{3/2} (k_B T)^{-(\delta_i+3)/2}
I^{\delta_i/2-1}\exp \left( -\frac{m_{i}|v-u|^{2}}{2k_{B}T}-\frac{%
I}{k_{B}T}\right).
\end{equation}
\end{enumerate}
\end{proposition}

\subsubsection{Mixture of polyatomic gases with discrete internal energies}\strut\\
We now extend the model described in Section \ref{sec:discrete-single} to a mixture of $N$ polyatomic gases \cite{Gio,Groppi-Spiga}, with molecular masses $m_1$, ..., $m_N$, where the polyatomicity of each species $i\in \{ 1,...,N\}$ is modeled by $N_{\intern,i}\ge 1$ different internal energies $I_i^{(1)}$, ..., $I_i^{(N_{\intern,i})} \in \RR_+$.

Observe that if $N_{\intern,1} = \dots= N_{\intern,N}=1$ (with $\varphi_{1}^{(1)}= \dots = \varphi_{N}^{(1)}=1$), the model reduces to the case of a mixture of monatomic species. We obtain a mixture model of monatomic and polyatomic species as soon as there exist at least one index $i$ such that $N_{\intern,i} = 1$, with $\varphi_{i}^{(k)}=1$, meaning that species $i$ is monatomic, and at least one index $j$ such that $N_{\intern,j}\ge2$, meaning that species $j$ is polyatomic.

Each collision can then be represented by two pre-collisional pairs and two corresponding post-collisional pairs, respectively indexed by $k$, $\ell$ and $k'$, $\ell'$, with $k$,  $k'\in \{1,\dots,N_{\intern,i}\}$, and $\ell$, $\ell'\in  \{1,\dots,N_{\intern,j}\}$.
More precisely, consider two colliding molecules of mass $m_i$ and $m_j$ with velocities and internal energies $( v,I_i^{(k)}) $
and $( v_*,I_j^{(\ell)}) $, changing into $( v',I_i^{(k')}) $ and $( v'_*,I_j^{(\ell')}) $, and define the internal energy gap
\begin{equation}\label{energy gap moment mixt discr int energy}
\Delta I_{ij}^{(k\ell,k'\ell')} =I_i^{(k')}+I_j^{(\ell')}-I_i^{(k)}-I_j^{(\ell)}.
\end{equation}%
The  microscopic momentum and energy conservations are given by
\begin{align}
m_i v+m_j v_* &= m_iv' +m_jv'_*, \phantom{\frac{m_i}2} \label{moment cons mixt discr int energy}\\
\frac{m_i}{2} |v|^2+\frac{m_j}{2}|v_*|^2 +I_i^{(k)} +I_j^{(\ell)} &=\frac{m_i}{2}|v'|^2+\frac{m_j}{2}| v'_*|^2+I_i^{(k')}+I_j^{(\ell')}.
\label{energy cons mixt discr int energy}
\end{align}%
Equivalently, the energy conservation can be rewritten in the center-of-mass reference frame, with the total energy 
\begin{equation*}
E_{ij}^{(k\ell)}=\dfrac{\mu _{ij}}{2} | V|^2+I_i^{(k)}+I_j^{(\ell)} =\dfrac{\mu _{ij}}{2} | V'|^2+I_i^{(k')}+I_j^{(\ell')}.
\end{equation*}
The post-collisional velocities are then given by
\begin{align*}
v'&=\dfrac{m_i v+m_j v_*}{m_i+m_j} +\sigma \dfrac{\mu_{ij}}{m_i}\sqrt{|V|^2-\dfrac{2}{\mu _{ij}}\Delta I_{ij}^{(k\ell,k'\ell')}}, \\ 
v'_*&=\dfrac{m_iv +m_j v_*}{m_i+m_j} -\sigma \dfrac{\mu_{ij}}{m_j}\sqrt{|V|^2-\dfrac{2}{\mu _{ij}}\Delta I_{ij}^{(k\ell,k'\ell')}}.
\end{align*}
The symmetry and microreversibility assumptions on the collision kernels for $( i,j) \in \{1,...,N\}^2$ are
\begin{gather*}
B_{ij}(v,v_*,I_i^{(k)},I_j^{(\ell)},I_i^{(k')},I_j^{(\ell')},\sigma)  
= B_{ji}( v_*,v,I_j^{(\ell)},I_i^{(k)},I_j^{(\ell')},I_i^{(k')},-\sigma) \\
=
B_{ij}\left( v',v'_*,I_i^{(k')},I_j^{(\ell')},I_i^{(k)},I_j^{(\ell)},\sigma'\right)
\\
B_{ii}( v,v_*,I_i^{(k)},I_i^{(\ell)},I_i^{(k')},I_i^{(\ell')},\sigma) 
= B_{ii}( v_*,v,I_i^{(\ell)},I_i^{(k)},I_i^{(k')},I_i^{(\ell')},\sigma).
\end{gather*}
The corresponding collision operator is written, for almost every $v\in\RR^3$, as
\begin{multline}
Q_{ij}^{(k)}(f,g)(v)=\sum_{k'=1}^{N_{\intern,i}}\sum_{\ell,\ell'=1}^{N_{\intern,j}}\int_{\RR^3\times\SS^2} \left( f^{(k')}(v') g^{(\ell')}(v'_*) \frac{\varphi_i^{(k)}\varphi_j^{(\ell)}}{\varphi_i^{(k')}\varphi _j^{(\ell')}} - f^{(k)}(v) g^{(\ell)}(v_*)\right) \\ \phantom{\int}
\times B_{ij}( v,v_*,I_i^{(k)},I_j^{(\ell)},I_i^{(k')},I_j^{(\ell')},\sigma)\varphi_i^{(k')}\varphi_j^{(\ell')}\frac{|V'|}{\left(E_{ij}^{(k\ell)}\right)^{1/2}} \,\md v_* \,\md\sigma,
\label{collisionoperator_mixture_discrete}
\end{multline}%
with constant positive degeneracies $\varphi_i^{(1)}$, ..., $\varphi_i^{(N_{\intern,i})}$ for each species $i$. 
In the previous equality \eqref{collisionoperator_mixture_discrete}, $f$ is the vector function consisting of the components for all different internal energies of species $i$, and $g$ the one consisting of the components for all different internal energies of species $j$.
In order to define the $i$th component of the vector collision operator, its argument has to be all these vector functions related to any species at the same time, as in the case of continuous energies \eqref{e:ithcollop}. Therefore, we use  the vector form of all components for different internal energies of all species 
$$f=\left( f_1^{(1)},\dots,f_1^{(N_{\intern,1})},f_2^{(1)},\dots,f_{N-1}^{(N_{\intern,N-1})}, f_N^{(1)},\dots,f_N^{(N_{\intern,N})}\right) ,$$
and using the same notation for a vector function $g$,
the $i$th component of the collision operator is given, for any $i\in\{1,\dots,N\}$, by
$$Q_{i}^{(k)} (f,g) = \sum_{j=1}^{N} Q_{ij}^{(k)} (f_i,g_j).$$
We can then set, for any $i$, 
\begin{equation} \label{e:defQimixtdiscr}
    Q_i(f,g)=\left(Q_{i}^{(1)}(f,g),\dots,Q_{i}^{(N_{\intern,i})}(f,g) \right),
\end{equation}
which allows to state a part of the $H$-theorem that defines the equilibrium state.

\begin{proposition}
\label{prop:htheorem_MDI} The three following properties are equivalent:
\begin{enumerate}
\item $\displaystyle Q_i(M,M)=0$ for any $i\in \{1,...,N\}$,

\item $\displaystyle\sum\limits_{i=1}^{N}\sum\limits_{k=1}^{N_{\intern,i}}\int_{\RR^{3}}Q_{i}^{(k)}(M,M)(v)\,%
\log \left[ \frac{M_{i}^{(k)}(v)}{ \varphi_{i} ^{(k)}} \right]\,\md v =0$,

\item there exist $n=(n_{1},\dots,n_N)\in \RR_+^N$, $u\in \RR^{3}$ and $T>0$
such that, for any $i\in \{1,...,N\}$ and $k\in \left\{
1,...,N_{\intern,i}\right\}$, and almost every $v$,
\begin{equation}
M_{i}^{(k)}(v)=\frac{n_{i}}{q_{i}}\left( \frac{m_{i}}{2\pi k_{B}T}\right)
^{3/2}\varphi _{i}^{(k)}\exp \left( -\frac{m_{i}|v-u|^{2}}{2k_{B}T}-\frac{%
I_{i}^{(k)}}{k_{B}T}\right) ,  \label{e:defMaxwMDI}
\end{equation}
with $q_{i}=\sum\limits_{k=1}^{N_{\intern,i}}\varphi _{i}^{(k)}\exp \left( -%
\dfrac{I_{i}^{(k)}}{k_{B}T}\right) $.
\end{enumerate}
\end{proposition}

\subsection{Boltzmann equation and linearized collision operator}
\begingroup 
The Boltzmann equation describes the time evolution of a system composed by a large number of particles, described by a distribution function $f$ defined on the phase space of the system. 
For the sake of simplicity, we also assume that the system is isolated, so that there is no external force acting on the particles.
The time evolution of the distribution function $f$ is governed by the Boltzmann equation
\begin{equation} \label{Boltz equ}
\frac{\partial f}{\partial t}+\left( v\cdot \nabla _x\right) f=Q\left( f,f\right),
\end{equation}
where $f$ denotes 
\setlength{\leftmargini}{1.5em}
\begin{compactitem}
\item the scalar function $f=f(t,x,v,I)$ in the case of a single polyatomic gas with continuous internal energy,
\item the vector function $f=\big(f^{(1)},\dots,f^{(N_\intern)}\big)$ in the discrete internal energy case, with $f^{(k)}=f^{(k)}(t,x,v)$,
\item the vector function  $f=\big(f_1, \dots, f_N \big)$, with $f_i=f_i(t,x,v)$ if $i\in \mathcal{M}$, and $f_i=f_i(t,x,v,I)$ if $i\in \mathcal{P}$, in the case of a mixture with continuous internal energies, 
\item the vector function $f=\big(f_1,\dots,f_N\big)$, where each $f_i$ is given by the vector function $f_i(t,x,v)=\big(f_i^{(1)}(t,x,v),\dots, f_i^{(N_{\intern,i})}(t,x,v)\big)$, in the discrete internal energy case.
\end{compactitem}

\smallskip

In each case, the collision operator $Q$ is a quadratic bilinear operator that accounts for the change of velocities and internal energies of particles due to the binary collisions.
\begin{compactitem}
\item For a single species, $Q$ is given by \eqref{collisionoperator} in the Borgnakke-Larsen framework, by \eqref{e:defQreson} in the resonant case, and by the vector form $Q = \big(Q^{(1)},\dots,Q^{(N_\intern)}\big)$, where each $Q^{(k)}$, $1\le k\le N_\intern$, is given by \eqref{collisionoperator_discrete} in the discrete internal energy case.
\item In the mixture case, when the internal energy variable is continuous, $Q$ is naturally defined as the vector expression  $Q=\big(Q_1,\dots,Q_N\big)$, where each $Q_i$ is given by \eqref{e:ithcollop}. 
\item In the mixture case, when the internal energy variable is discrete, $Q$ is again defined as $Q=\big(Q_1,\dots,Q_N\big)$, where, this time, each $Q_i$ is given by \eqref{e:defQimixtdiscr}.
\end{compactitem}

\smallskip 

In order for us to describe the linearization setting in a unified way, let us rewrite the scalar collision operator in the following way, for almost every $w=(v,I)$, 
\begin{equation} \label{e:Qgeneral}
    Q(f,g) (w) = \int_D \Big( f' g'_* \Phi - f g_*\Big) A(w,W) \,\md W,
\end{equation}
where the integration variable $W$ is
\begin{equation}\label{e:W}
W= \begin{cases}
  (v_*,I_*,I',\sigma) & \text{for the resonant case,}\\
  (v_*,I_*,R,r,\sigma) & \text{for the Borgnakke-Larsen framework,}\\
    \end{cases}
\end{equation}
the integration domain $D$ is either $D=\RR^3\times (\RR_+)^2\times \SS^2$ or $D=\RR^3\times (\RR_+)\times(0,1)^2\times \SS^2$, $\Phi= \left[(II_*)/(I'I'_*)\right]^{\delta/2-1}$, which is obviously similar for all continuous internal energy models, and, up to some multiplicative constants if necessary,
\begin{equation}\label{e:A}
A(w,W)= \begin{cases}
  B(v,v_*,I,I_*,I',\sigma)  \,\frac{\left[ I'(I+I_*-I') \right]^{\delta/2-1}}{(I+I_*)^{\delta-1}} & \text{(resonant),}\\
  B(v,v_*,I,I_*,r,R,\sigma) \\\qquad \times\,\left( r (1-r) \right)^{\delta/2 - 1} (1-R)^{\delta - 1} \sqrt{R} & \text{(Borgnakke-Larsen).}\\
    \end{cases}
\end{equation}

\medskip

In the case of vector quantities, \eqref{e:Qgeneral} can be extended in a straightforward way. First, in the discrete internal energy case, $Q^{(k)}$ writes, for almost every $w$, 
\begin{multline}\label{e:Qgeneral_discr}
    Q^{(k)} (f,g) (w) \\= \sum_{k',\ell,\ell'=1}^N  \int_D \Big( f^{(k')}(v') g^{(\ell')}(v'_*) \Phi^{(k\ell,k'\ell')} - f^{(k)}(v) g^{(\ell)}(v_*)\Big) A^{(k\ell,k'\ell')} (w,W)\,\md W,
\end{multline}
where $w=v$, 
$W=(v_*,\sigma)$, $D = \RR^3 \times \SS^2$, $\Phi^{(k\ell,k'\ell')} = (\varphi^{(k)} \varphi^{(\ell)})/( \varphi^{(k')} \varphi^{(\ell')} )$ and
\begin{equation}\label{e:A-discrete}
A^{(k\ell,k'\ell')} (w,W)=B( v,v_*,I^{(k)},I^{(\ell)},I^{(k')},I^{(\ell')},\sigma )  \frac{\varphi^{(k')} \varphi^{(\ell')} |V'|}{(E^{(k\ell)})^{1/2}}.
\end{equation}
Second, in the mixture case with a continuous internal energy variable, $Q_i$ writes, for almost every $w$, 
\begin{equation}\label{e:Qgeneral_mixt}
    Q_{i}(f,g)(w) = \sum_{j=1}^N \int_D \Big(f'g'_* \Phi_{ij} - fg_*\Big) A_{ij}(w,W) \,\md W,
\end{equation}
where $w = (v,I)$ if $i\in\mathcal{P}$ and $w=v$ if $i\in\mathcal{M}$, $W$ is defined by 
\begin{equation}\label{e:W-mix}
    W = \begin{cases}
    (v_*,I_*,R,r,\sigma) & \text{if $i,j\in \mathcal{P}$,}\\
    (v_*,R,\sigma) & \text{if $i\in \mathcal{P}$, $j\in \mathcal{M}$,}\\
    (v_*,I_*,R,\sigma) & \text{if $i \in \mathcal{M}$, $j\in \mathcal{P}$,}\\
    (v_*,\sigma) & \text{if $i,j\in \mathcal{M}$,}\\
    \end{cases}
\end{equation}
with the relevant corresponding domain $D$, $\Phi_{ij}$ is defined by
\begin{equation*}
    \Phi_{ij} = \begin{cases}
    \left(\frac{I}{I'}\right)^{\delta_i/2-1} \left(\frac{I_*}{I'_*}\right)^{\delta_j/2-1} & \text{if $i,j\in\mathcal{P}$,} \\
    \left(\frac{I}{I'}\right)^{\delta_i/2-1}  & \text{if $i\in\mathcal{P}$, $j\in\mathcal{M}$,} \\
     \left(\frac{I_*}{I'_*}\right)^{\delta_j/2-1} & \text{if $i\in\mathcal{M}$, $j\in\mathcal{P}$,} \\
    1 & \text{if $i,j\in\mathcal{M}$,}
    \end{cases}
\end{equation*}
and $A_{ij}(w,W)$ by
\begin{equation}\label{e:Aij}
    A_{ij}(w,W) = \begin{cases}
    B_{ij}(v, v_*, I, I_*, r, R, \sigma) \\\qquad \times\, r^{\delta_i/2-1}(1-r)^{\delta_j/2-1} \, (1-R)^{\delta_i/2 + \delta_j/2-1} \, \sqrt{R}, & \text{if $i,j\in \mathcal{P}$,}\\
    B_{ij}(v, v_*, I, R, \sigma) \,  (1-R)^{\delta_i/2 -1} \, \sqrt{R}, & \text{if $i\in \mathcal{P}$, $j\in \mathcal{M}$,}\\
    B_{ij}(v, v_*, I_*, R, \sigma) \,  (1-R)^{\delta_j/2-1} \, \sqrt{R}, & \text{if $i \in \mathcal{M}$, $j\in \mathcal{P}$,}\\
    B_{ij}(v, v_*,\sigma)  & \text{if $i,j\in \mathcal{M}$.}\\
    \end{cases}
\end{equation}

Last, the case of a mixture of polyatomic gases with discrete internal energies is a straightforward extension of \eqref{e:Qgeneral_discr} and \eqref{e:Qgeneral_mixt} by writing, for almost every $w$,
\begin{multline*}
Q_i^{(k)}(f,g)(w) \\ = \sum_{j=1}^N \sum_{k'=1}^{N_{\intern,i}}\sum_{\ell,\ell'=1}^{N_{\intern,j}} \int_D  \Big( f_i^{(k')}(v') g_j^{(\ell')}(v'_*) \Phi_{ij}^{(k\ell,k'\ell')} - f_i^{(k)}(v) g_j^{(\ell)}(v_*)\Big) A_{ij}^{(k\ell,k'\ell')} (w,W)\,\md W, 
\end{multline*}
where $w=v$, $W=(v_*,\sigma)$, $D=\RR^3\times\SS^2$, $\Phi_{ij}^{(k\ell,k'\ell')} = (\varphi_i^{(k)} \varphi_j^{(\ell)})/( \varphi_i^{(k')} \varphi_j^{(\ell')})$ and 
\begin{equation}\label{e:Aij-discrete}
    A_{ij}^{(k\ell,k'\ell')} (w,W) =B_{ij}( v,v_*,I_i^{(k)},I_j^{(\ell)},I_i^{(k')},I_j^{(\ell')},\sigma )  \frac{\varphi_i^{(k')} \varphi_j^{(\ell')} |V'|}{(E_{ij}^{(k\ell)})^{1/2}}.
\end{equation}
\smallskip

Further, the associated equilibria for a single species are Maxwellian distributions $M$ given by \eqref{e:defMaxwBL} in the Borgnakke-Larsen framework, \eqref{e:defMaxwresonpowerlaw} in the resonant case, and \eqref{e:defMaxwDI} when handling discrete internal energy levels.
For mixtures, the Maxwellian distributions $M_i$ are defined in \eqref{e:defMaxwMBLmono}--\eqref{e:defMaxwMBLpoly} for the Borgnakke-Larsen model, and \eqref{e:defMaxwMDI} in the discrete internal energy case.
\endgroup

The standard perturbative setting for the Boltzmann equation \eqref{Boltz equ}, combined with the $H$-theorem, leads to considering deviations of Maxwellian distributions under the form
\begin{equation}
f=M+M^{1/2}h.  \label{linear}
\end{equation}
In this case, the linearized Boltzmann operator is defined as
\begin{equation} \label{e:L}
    \mathcal L h = M^{-1/2} \left[ Q(M,M^{1/2} h) + Q (M^{1/2} h,M)\right].
\end{equation}
In the case of mixtures, this equation has to be understood in the following sense for any $1\leq i \leq N$
\begin{equation}\label{e:L_mixt}
[\mathcal L h]_i = M_i^{-1/2} \sum_{j=1}^N \left[ Q_{ij}(M_i,M_j^{1/2} h_j) + Q_{ij} (M_i^{1/2} h_i,M_j)\right].
\end{equation}
This linearized Boltzmann operator can be written as $\mathcal L =  K - \nu \text{Id}$, where the collision frequency $\nu$ is defined with the notations of \eqref{e:Qgeneral} as
\begin{equation*}
    \nu (w) = \int_D M_* A(w,W)\,\md W. 
\end{equation*}
The extension to vector quantities, with the notations of \eqref{e:Qgeneral_discr} or \eqref{e:Qgeneral_mixt}, is straightforward.

The operator $K$ consists of three contributions and is given with the notations of \eqref{e:Qgeneral} by
\begin{equation}\label{e:opK}
    Kh (w) = M^{1/2} \int_D \Big( (M^{-1/2} h)' + (M^{-1/2} h)'_* - (M^{-1/2} h)_*\Big) M_* A(w,W) \,\md W. 
\end{equation}
The extension to discrete internal energies is similar as in \eqref{e:Qgeneral_discr}, whereas the extension to mixtures has to be understood as in \eqref{e:Qgeneral_mixt} and \eqref{e:L_mixt}.

\section{Compactness property of \textit{K}} \label{s:compactness}

Recent works have tackled the question of the compactness of the operator $K$, defined in \eqref{e:opK}, for the different polyatomic models presented in Section \ref{sec:models}. In these papers, the compactness property is proved under some assumptions on the collision kernels, basically imposing some bounds on their growth.

\subsection{Monatomic gases}
Let us briefly recall under which assumption the operator $K$ is compact for monatomic gases, in the case of hard potential and with Grad's cut-off assumption, since it is the framework of all papers  which tackled this question for polyatomic gases. For more references on the monatomic case, see the review paper \cite{BS-16}. For both monatomic single species and mixtures of monatomic gases, compactness of the operator
~$K$ has been proved \cite{Gr-63,BGPS-13} under the assumption \hypref{hyp:mono} stated below.
\begin{hyp}\label{hyp:mono}
There exist a constant $C>0$ and an exponent $0<\zeta <1$ such that the possible collision kernels $B$ or $B_{ij}$, for any $1\leq i,j \leq N$, satisfy
\begin{equation*} 
B(v,v_*,\sigma),~ B_{ij}(v,v_*,\sigma) \leq C\left\vert V\right\vert \left( 1+\dfrac{1}{\left\vert V\right\vert^{2-\zeta }}\right).
\end{equation*}
\end{hyp}
\noindent Since monatomic gases are discussed in detail in \cite{BS-16}, we do not discuss them here.
\medskip 

For a single polyatomic gas, we described three models: one written within the Borgnak\-ke-Larsen framework, one for resonant collisions, and one with discrete internal energies. In those three cases, the assumptions on the collision kernels made by the authors are slightly different. Let us state them below.

\subsection{Single polyatomic gas with continuous internal energy}
For the model written with the Borgnakke-Larsen procedure, several contributions appeared recently, obtained under different assumptions. 

\begin{hyp}\label{hyp:niclas}
There exist a constant $C>0$ and an exponent $0<\zeta <1$ such that
\begin{equation*}
B(v,v_*,I,I_*,r,R,\sigma)\leq CE\left( 1+\dfrac{1}{\left( \left\vert V\right\vert \left\vert
V^{\prime }\right\vert \right) ^{1-\zeta /2}}\right).
\end{equation*}
\end{hyp}
\begin{hyp}\label{hyp:brull}
Let $\zeta>-1$ and $\delta\geq 2$.
There exist a constant $C>0$ and a function $\Psi(r,R)$, satisfying the symmetry condition $\Psi \left( r,R\right) =\Psi \left( 1-r,R\right) $ and the integrability condition
\begin{equation*}
    \Psi\left( r,R\right)^{2} \left( 1-r\right) ^{\delta -3-\zeta }r^{\delta/2-2}R(1-R)^{3\delta /2-3-\zeta }  \in L^{1}\left( \left( 0,1\right)^{2}\right) ,
\end{equation*}
such that 
\begin{equation*}
    B(v,v_*,I,I_*,r,R,\sigma) \leq C\Psi \left( r,R\right) E^{\zeta /2}.
\end{equation*}
\end{hyp}
In this case, compactness of $K$ has been proved in two different contributions.
\begin{theorem}[proved in \cite{Be-23b}]\label{th:niclas}
For any $\delta\geq 2$, under Hypothesis \hypref{hyp:niclas}, the operator $K$ is compact from $L^2(\RR^3\times \RR_+)$ into itself.
\end{theorem}

\begin{theorem}[proved in \cite{BST-24a,shahine-phd}] \label{th:brull}
For any $\zeta > -1$, for $\delta\geq 2$ and $B$ satisfying Hypothesis \hypref{hyp:brull}, the operator $K$ is compact from $L^2(\RR^3\times \RR_+)$ into itself.
\end{theorem}
Observe that Hypothesis \hypref{hyp:brull} couples the value of $\delta$ with the assumption on the collision kernel $B$, through the integrability of the function $\Psi$. Let us also mention that in the diatomic case ($\delta=2$), the operator $K$ was proved to be compact in \cite{BST-22}, but under more restrictive assumptions than Hypothesis \hypref{hyp:brull} for $\delta=2$.

Let us also underline that a discussion on these different assumptions is led in Section~\ref{sec:discussion}. In particular, for a comparison of \hypref{hyp:niclas} and \hypref{hyp:brull}, we thus refer to this section.

\begin{remark}
By Theorem \ref{th:niclas-mix} restricted to a single species, Theorem \ref{th:niclas} remains valid for $E$ replaced by $E^{1/2}$ in Hypothesis \hypref{hyp:niclas}. This extends Theorem \ref{th:niclas} to be valid not only for hard and super-hard potential like kernels, but also for some soft potential like kernels.
\end{remark}

\subsection{Single polyatomic gas with resonant collisions}

In the resonant case, the collision kernel is assumed \cite{borsoni-phd,BBS-23} to be upper-bounded by a tensored form on $
(|V|,\cos \theta )$ on the one hand and on $(I,I_{\ast },I^{\prime })$ on the other hand. 
\begin{hyp}\label{hyp:resonant}
There exist two functions $b_\kin$ and $b_\intern$ satisfying, for some constant $C>0$ and some exponents $\zeta \in \lbrack 0,1)$, $\zeta_{1}\in \lbrack 0,1/2)$ and $\zeta_{2}\in (-\delta,\delta)$, 
\begin{align}
0 \leq b_\kin(|V|,\cos \theta )& \leq C\left[ |\sin \theta |\left(
|V|^{2}+|V|^{-1}\right) +|V|+|V|^{-\zeta }+|\sin \theta |^{-\zeta _{1}}%
\right] \phantom{\frac12} \nonumber \\
0 \leq b_\intern(I,I_{\ast })& \leq C(I+I_{\ast })^{1+\zeta _2/2-\delta}, \label{e:hypbint}
\end{align}%
such that
\begin{equation} \label{e:tensformBres}
B(v,v_*,I,I_*,I',\sigma)\leq b_\kin(|V|,\cos \theta )\, b_\intern(I,I_{\ast })\, \mathbf{1}_{[0,I+I_*]}(I'),
\end{equation}
or any linear combination of such terms.
\end{hyp}
Under this assumption, the compactness of $K$ can be proved.
\begin{theorem}[proved in \cite{BBS-23}]\label{th:resonant}
Under Hypothesis \hypref{hyp:resonant}, the operator $K$ is compact from $L^2(\RR^3\times \RR_+)$ into itself.
\end{theorem}

\subsection{Single polyatomic gas with discrete internal energy}
In the case of discrete internal energies, an analogous result to Theorem \ref{th:niclas} has been proved under the following assumption.
\begin{hyp}\label{hyp:discrete}
There exist a constant $C>0$ and an exponent $0<\zeta <1$ such that, for any $1 \leq k,k',\ell,\ell' \leq N_\intern$,
\begin{equation*}
B(v,v_*,I^{(k)},I^{(\ell)},I^{(k')},I^{(\ell')},\sigma)\leq C\left( E^{\left(k\ell\right)} \right)^{1/2} \left( 1+\dfrac{1}{\left(
\left\vert V\right\vert \left\vert V'\right\vert \right) ^{1-\zeta
/2}}\right) .
\end{equation*}
\end{hyp}
\begin{theorem}[proved in \cite{Be-23a}]
Under Hypothesis \hypref{hyp:discrete}, the operator $K$ is compact from $L^2(\RR^3)^{N_\intern}$ into itself.
\end{theorem}

\subsection{Mixture of monatomic and polyatomic gases with continuous internal energy}

When considering a mixture with polyatomic gases, several configurations can occur, if also monatomic gases exist in the mixture. Moreover, as for the case of a single species, different contributions appeared in the literature, using different assumptions. These contributions are natural extensions of the results we mentioned in the previous subsection about the single species case.
\begin{hyp}\label{hyp:niclas-mix}
There exist a constant $C>0$ and an exponent $0<\zeta <1$ such that the possible $B_{ij}$, for any $1\leq i,j\leq N$, satisfy
\begin{multline*}
B_{ij}(v,v_*,I,I_*,r,R,\sigma), ~ B_{ij}(v,v_*,I,R,\sigma),~ B_{ij}(v,v_*,I_*,R,\sigma),~ B_{ij}(v,v_*,\sigma)\\
\leq CE_{ij}^{1/2}\left( 1+\dfrac{1}{\left( \left\vert V\right\vert
\left\vert V^{\prime }\right\vert \right) ^{1-\zeta /2}}\right).
\end{multline*}
\end{hyp}

\begin{hyp}\label{hyp:brull-mix}
Let $\zeta_{ij} > -1$ and $\delta_i\geq 2$, for any $1\leq i,j\leq N$, such that $\delta_{i}-\delta_{j}\leq 2+\zeta _{ij}$. There exist a constant $C>0$ and functions $\Psi_{ij}(r,R)$ satisfying the symmetry condition $\Psi _{ij}( r,R) =\Psi _{ij}( 1-r,R)$  and the integrability conditions
\begin{gather*}
\Psi _{ij}( r,R)^{2} ( 1-r) ^{\delta_{j}/2-2}r^{( \delta _{i}+\delta_{j}) /2-3-\zeta_{ij}}R(1-R)^{\delta _{i}/2+\delta _{j}-3-\zeta _{ij}} \in L^{1}( (0,1)^{2}),  \\
\Psi _{ij}( r,R)^{2} ( 1-r) ^{\delta _{j}-3-\zeta_{ij}}r^{\delta_{i}/2-2}R(1-R)^{\delta _{i}/2+\delta _{j}-3-\zeta _{ij}} \in L^{1}( (0,1)^{2}) ,
\end{gather*}
such that
\[B_{ij}(v,v_*,I,I_*,r,R,\sigma) \leq C\,\Psi _{ij}( r,R)\, E_{ij}^{\zeta _{ij}/2}.\]
\end{hyp}
The additional assumption  $\delta_{i}-\delta_{j}\leq 2+\zeta _{ij}$ is not explicitly stated in \cite{shahine-phd}, but the current proof seems to impose this condition.

\begin{theorem}[proved in \cite{Be-24a}]\label{th:niclas-mix}
For any $(\delta_i)_{1\leq i \leq N}$ with $\delta_i \geq 2$, under Hypothesis \hypref{hyp:niclas-mix}, the operator $K$ for a mixture of polyatomic gases involving possibly also monatomic gases is compact from $L^2(\RR^3)^{|\mathcal{M}|} \times L^2(\RR^3\times \RR_+)^{|\mathcal{P}|}$ into itself, up to an index reordering.
\end{theorem}

\begin{theorem}[proved in \cite{BST-24c}, \cite{shahine-phd}]\label{th:brull-mix}
For any $(\zeta_{ij})_{1\leq i,j\leq N}$ with $\zeta_{ij}>-1$, for $(\delta_{i})_{1\leq i\leq N}$ and $B_{ij}$ satisfying Hypothesis \hypref{hyp:brull-mix}, the operator $K$ for a mixture of polyatomic gases is compact from $L^2(\RR^3\times \RR_+)^N$ into itself.
\end{theorem}
Observe that, as in Hypothesis \hypref{hyp:brull}, \hypref{hyp:brull-mix} couples the values of $\delta_{i}$ and $\delta_j$ with the assumption on the collision kernel $B_{ij}$, through the integrability of the function $\Psi_{ij}$.

\begin{remark}
For an extension of Theorem \ref{th:niclas-mix} to a model including chemical reactions in form of dissociation and association (with possible extension to other chemical reactions through instant combinations of such), we refer to a recent paper \cite{Be-24c}.
\end{remark}

\subsection{Mixture of monatomic and polyatomic gases with discrete internal energy}
As before, in the discrete internal energy case, the single species case is also extended to the mixtures of polyatomic gases, with possibly monatomic species. The assumption on the collision kernel is adapted as follows.
\begin{hyp}\label{hyp:discrete-mix}
There exist a constant $C>0$ and an exponent $0<\zeta <1$ such that for any $1 \leq k,k'\leq N_{\intern,i}$, $1 \leq \ell,\ell'\leq N_{\intern,j}$, 
\begin{equation*}
B_{ij}(v,v_*,I_i^{(k)},I_j^{(\ell)},I_i^{(k')},I_j^{(\ell')},\sigma)\leq C\left( E_{ij}^{\left(k\ell\right)} \right)^{1/2} \left( 1+\dfrac{1}{\left(\left\vert V\right\vert \left\vert V'\right\vert \right) ^{1-\zeta/2}}\right).
\end{equation*}
\end{hyp}
The compactness result follows. 
\begin{theorem}[proved in \cite{Be-24b}]
Under Hypothesis \hypref{hyp:discrete-mix}, the operator $K$ for a mixture of polyatomic gases involving possibly also monatomic gases, is compact from $L^2(\RR^3)^{N_\intern^\textnormal{tot}}$, where $N_\intern^\textnormal{tot} = \sum_{i=1}^N N_{\intern,i}$, into itself.
\end{theorem}

\section{Main ideas for proving compactness in the single species case}\label{s: main ideas}

From the general form \eqref{e:opK} of the operator $K$, we observe that there are three different contributions
\begin{align*}
    K_1 h(w) &= - M(w)^{1/2} \int_D  h(w_*) M(w_*)^{1/2} A(w,W) \,\md W,\\
    K_2 h(w) &= M(w)^{1/2} \int_D M(w'_*)^{-1/2} h(w'_*) M(w_*) A(w,W) \,\md W,\\
    K_3 h(w) &= M(w)^{1/2} \int_D M(w')^{-1/2} h(w') M(w_*) A(w,W) \,\md W.
\end{align*}

The mutual strategy used in the different approaches is to find a kernel form of each operator and prove some bounds and integrability properties of the kernel, which allow to deduce compactness for the corresponding operator.

The first observation is that the operator $K_1$ can be treated in a straightforward way, since it is already under a kernel form, remembering that the variable $W$ contains $w_*$ and parameters. 
The integrability of this kernel is easy to obtain, under classical assumptions on the collision kernel, and is a natural extension of the monatomic case \cite{Gr-63,BS-16}.
Indeed, noting for example that for a polyatomic single species (a multiple of) the exponent of the product $MM_{\ast }$ can be recast as%
\begin{equation*}
m\frac{\left\vert v\right\vert ^{2}}{2}+m\frac{\left\vert v_{\ast
}\right\vert ^{2}}{2}+I+I_{\ast }=m\left\vert \frac{v+v_{\ast }}{2}%
\right\vert ^{2}+E,
\end{equation*}%
combined with estimates of the form
\begin{equation}
I,I_{\ast },I^{\prime },I_{\ast }^{\prime }\leq C E\text{,}  \label{ineq}
\end{equation}%
the kernel, by a suitable change of variables, can be proved to be in $L^{2}((\RR^3\times \RR_+)^2)$.
Hence, as a Hilbert-Schmidt integral operator, $K_1$ is compact.

The main contribution is thus to treat the two remaining operators $K_2$ and $K_3$. In the single species case, it is worth noticing that both operators have the same structure, and thus similar proofs can be used to prove their compactness. We will thus only detail the treatment of $K_2$ or $K_3$ in this section.

In the following, we will explain the main ideas of the different contributions on the compactness of the operator $K$. To this end, we will drop any multiplicative constant which do not play any role in the proof. 
In particular, without loss of generality, the Maxwellian functions can be assumed to be centered and normalized, giving some additional constant in front of the integral.
The notations $\propto$ and $\lesssim$ will be used to denote an equality and an inequality up to (positive) multiplicative constants, which may depend on the various fixed parameters.

\subsection{Sketch of the proof of Theorem~\ref{th:brull}} \label{ss:brullmonospecies}
Let us describe the main ideas of Brull, Shahine and Thieullen's approach \cite{BST-22,BST-24a
}. 
We recall the expression of $K_2$, using the expressions \eqref{e:A} of $A(w,W)$ and \eqref{e:defMaxwBL} of $M$ together with the relation $M M_* = M' M'_* \Phi$
\begin{multline}\label{e:opK'* Brull}
    K_2h (v,I) \propto I^{\delta/4-1/2} \int_D h'_* e^{-\tfrac{1}{4}(|v'|^2+|v_*|^2)} e^{-\tfrac{1}{2}(I'+I_*)} (I'_*)^{-\delta/4+1/2} I_*^{\delta/2-1} 
    \\
\times    (r (1-r))^{\delta/2-1} (1-R)^{\delta-1} \sqrt{R} B   \,\md v_* \, \md I_* \, \md r \, \md R \, \md \sigma. 
\end{multline}

In order to get a kernel form, it is natural to change the variables $(v_*,I_*)$ into $(v'_*,I'_*)$, for which the Jacobian $J=8/(1-r)(1-R)$ can be computed explicitly. From the collision rules, explicit expressions of $v'$, $I'$, $v_*$, $I_*$ depending on $v$, $I$, $v'_*$, $I'_*$, $r$, $R$, $\sigma$ are obtained. The nonnegativity of internal energies leads to one nontrivial condition ($I_*\geq 0$), which can be translated into a restriction on the domain of integration $D$. This restriction is first written explicitly as one on the domain of $(v'_*,I'_*)$ only. A change of perspective on this restricted domain allows to write it (implicitly) as a restriction on the domain of the parameters $r,R,\sigma$, while the domain of $(v'_*,I'_*)$ is the original full one $\RR^3\times \RR_+$. 
In this way, the operator is written under a kernel form $\int_{\RR^3\times \RR_+} h(v'_*,I'_*) k_2(v,I,v'_*,I'_*)\,\md v'_* \,\md I'_*$, where $k_2$ contains an integral on the restricted domain for $r,R,\sigma$. 

In \cite{BST-24a}, the authors then prove that the kernel is integrable with respect to all variables, \ie $k_2 \in L^2((\RR^3\times \RR_+)^2)$.
First, the kernel is bounded by the same integral on the whole domain $(0,1)^2\times \SS^2$ for the parameters $r,R,\sigma$, and Cauchy-Schwarz is used in a standard way, giving the following estimate of the norm of $k_2$ in $L^2((\RR^3\times \RR_+)^2)$ 
\begin{multline} \label{e:k2L2norm}
    \|k_2\|_2^2 \lesssim  \int_{\RR^3\times \RR_+ \times D} {I^{\delta/2-1}} e^{-\tfrac{1}{2}(|v'|^2+|v_*|^2)} e^{-(I'+I_*)} (I'_*)^{-\delta/2+1} I_*^{\delta-2} 
    \\
\times    (r (1-r))^{\delta-2} (1-R)^{2\delta-2} R B^2 J^2   \,\md v \, \md I\,\md v'_* \, \md I'_* \, \md r \, \md R \, \md \sigma. 
\end{multline}
Using that $I'_* = (1-r)(1-R)E$, the term $(I'_*)^{-\delta/2+1}$ gives a factor $E^{-\delta/2+1}$ with the corresponding corrections on the powers of $r$ and $R$. Moreover, since $I \leq E$, this negative power of $E$ cancels with the very same power of $I$, provided $\delta/2-1\geq 0$ (see Hypothesis \hypref{hyp:brull}). Then, the inverse change of variables $(v'_*,I'_*) \mapsto (v_*,I_*)$
is used, as well as the assumption \hypref{hyp:brull} on $B$, which gives a contribution $E^\zeta$. Using the expression of $J$ and $I'$, the change of variable $I \mapsto E$ in \eqref{e:k2L2norm} leads to
\begin{multline*}
    \|k_2\|_2^2 \lesssim \int_{(0,1)^2}       
  r^{\delta-2} (1-r)^{\delta/2-2} (1-R)^{3\delta/2-2} R \Psi(r,R)^2   \\
\times \int_{\RR_+}E^{\zeta} e^{-r(1-R)E}
\int_{\RR^3\times \RR_+}I_*^{\delta-2} e^{-I_*} e^{-\tfrac{1}{2}|v_*|^2}
\left(\int_{\RR^3\times \SS^2} e^{-\tfrac{1}{2}|v'|^2} \,\md v \, \md \sigma\right)
\,\md v_* \, \md I_*\, \md E \, \md r \, \md R . 
\end{multline*}
For fixed $E,R,v_*$, the change of variables $v \mapsto v'$ allows to bound the integral in the parentheses. Further, the integral with respect to $v_*,I_*$ is obviously bounded. The integration in $E$ is possible for $\zeta>-1$ (\textit{cf.} Hypothesis \hypref{hyp:brull}) and gives an additional singularity in $r$ and $1-R$. The assumptions on $\Psi$ given in Hypothesis \hypref{hyp:brull} are finally chosen in order to have integrability of the rest with respect to $r,R$.

There is no additional difficulty to treat $K_3$, since the same ideas apply. The change of variables is obviously $(v_*,I_*) \mapsto (v',I')$ and some powers of $r,1-r,R,1-R$ are modified (due to the Jacobian of the change of variables and the expressions of $I'$ and $I'_*$).

\subsection{Sketch of the proof of Theorem \ref{th:niclas}}
Let us describe the main ideas of Bernhoff's approach \cite{Be-23b}.
The paper is written with a non-parametrized version of the collision operator. Nevertheless, with a suitable parametrization, the Borgnakke-Larsen operator considered here is obtained. Moreover, during the process of writing the integral operators $K_{1},K_{2},K_{3}$ under a kernel form, some changes of notation on the velocities are used to obtain a unified kernel form, \ie with the same notation for the independent variables of the kernel. However, those changes of notation are avoided here for the sake of consistency and clarity.

Concerning the operator $K_{3}$ (or, the operator $K_{2}$ by shifting the roles of the
variables $w'$ and $w'_{\ast }$), a new
parametrization is used
\begin{equation*}
\left\{ 
\begin{array}{c}
v_{\ast }=v^{\prime }+\xi-\chi n \\ 
v_{\ast }^{\prime }=v+\xi-\chi n%
\end{array}%
\right. , \quad \text{ with } \quad 
n=\frac{v-v'}{|v-v'|}, \qquad \xi\perp n ,
\end{equation*}%
and $\chi =\Delta I/ (m|v-v'|)$, where the energy gap is given by $ \Delta I = I'_*+I'-I_*-I$. This parametrization results in the following decomposition of $v+v'$
\begin{equation*}
\left( v+v^{\prime }\right) _{n}:=\left( v+v^{\prime }\right) \cdot n=\frac{%
\left\vert v\right\vert ^{2}-\left\vert v^{\prime }\right\vert ^{2}}{%
\left\vert v-v^{\prime }\right\vert }\quad \text{ and }\quad \left( v+
v^{\prime }\right) _{\perp _{n}}:=v+v^{\prime }-\left(
v+v^{\prime }\right) _{n}n\text{.}
\end{equation*}%
Note that for a vanishing internal energy gap, this is the parametrization
of Grad \cite{Gr-63,Glassey} for a monatomic gas. Then (a multiple of) the exponent of the
product $M_{\ast }M'_{\ast }$ can be recast as%
\begin{multline*}
\frac{\left\vert v_{\ast }\right\vert ^{2}}{2}+\frac{\left\vert v_{\ast
}^{\prime }\right\vert ^{2}}{2}+\frac{I_{\ast }}{m}+\frac{I_{\ast }^{\prime }%
}{m} \\
=\left\vert \frac{\left( v+v^{\prime }\right) _{\perp _{n}}}{2}+\xi\right\vert ^{2}\!+\frac{\left( \left\vert v^{\prime }\right\vert
^{2}-\left\vert v\right\vert ^{2}+2\chi \left\vert v-v^{\prime }\right\vert
\right) ^{2}}{4\left\vert v-v^{\prime }\right\vert ^{2}}+\frac{\left\vert
v-v^{\prime }\right\vert ^{2}}{4}+\frac{I_{\ast }}{m}+\frac{I_{\ast
}^{\prime }}{m}\text{.}
\end{multline*}%
Similar arguments as for a monatomic gas \cite{Gr-63,Glassey} can be applied for the velocity
part. For the internal energy part, the integration domains of $I$ and $I'$ are split into upper and lower intervals followed by convenient estimates of the form \eqref{ineq} being applied on the different domains. The internal energy factor of the upper bound of the kernel to be considered is of the form
\begin{equation*}
\widetilde{k}\left( I,I'\right) =\left( II'\right)^{\delta /4-1/2}\int_{0}^{+\infty }\int_{0}^{+\infty }e^{-\left( I_{\ast}+I'_{\ast }\right) /2}\frac{\left( I_{\ast }I'_{\ast }\right) ^{\delta /2-1}}{E^{\delta -1/2}}\,\md I_{\ast } \md I'_{\ast }.
\end{equation*}
Although $E$ is depending on velocity, it is estimated below by $I, I_*, I'$, or $I'_*$, and discussions on the exponents allow to prove that the desired integrability properties of $\widetilde{k}$.

To obtain the compactness of the
operator, it is proved to be a uniform limit of Hilbert-Schmidt integral
operators, based upon the kernel form of the operator and the observations mentioned before. The strategy is to
prove that 
\begin{enumerate}
    \item the integral of the kernel with respect to $w$ is bounded in $w'$; 
    \item the kernel is $L^{2}$ in an increasing sequence of truncated domains; 
    \item the supremum (over the domain of $w$) of the integral of the kernel over the complements of the truncated domains with respect to $w'$ is vanishing for the limiting sequence.
\end{enumerate} 

We recall that the comparison of the assumptions needed to prove compactness in the two approaches, described in Section~\ref{ss:brullmonospecies} and in this section, are discussed in Section~\ref{sec:discussion}. Nevertheless, at this point, we can wonder what differs between the proofs of the two approaches, and what could allow to have a wider possible range of parameters. 
One observation is that Bernhoff's approach is in a sense similar to the classical one for the monatomic linearized operator. Indeed, when the internal energies (and internal energy gaps) vanish, Grad's parametrization and scheme of the proof are recovered. 
Another observation is that in this approach, the operator is proved to be compact without the kernel necessarily being in $L^2((\RR^3\times \RR_+)^2)$. With the estimates of \cite{Be-23b}, when trying to prove that the kernel $k_2$ being in $L^2((\RR^3\times \RR_+)^2)$, some more restrictive assumptions are needed, in particular  $\delta$ having to be sufficiently large.

\subsection{Sketch of the proof of Theorem \ref{th:resonant}}

Let us describe the main ideas for the resonant case \cite{BBS-23}. In that situation, we recall that the microscopic velocities and internal energies evolve separately, and that Hypothesis~\hypref{hyp:resonant} on the collision kernel $B$ provides an upper-bound of $B$ as a product of $b_\kin$, only depending on the velocity variables and parameter, and $b_\intern$ multiplied by a characteristic function, both only involving internal energies. This means that one can deal with the resonant $K_2$ and $K_3$ by separately treating the velocity and internal energy parts of the operator. Let us be more accurate on $K_3$. 

The writing of $K_3$ into a kernel form is obtained in the same way as Grad proceeds, see \cite{Gr-63, BGPS-13}. Indeed, Grad's kernel derivation only involves velocity quantities, leading to a kernel $k_3(v,I,\eta,I')$. Thanks to the tensored upper-bound \eqref{e:tensformBres} of $B$, it is easy to obtain a bound of the kernel 
of the kind
\begin{equation*}\label{e:bornek3res}
    k_3(v,I,\eta,I') \lesssim k_\kin (v,\eta) k_\intern (I,I'), 
\end{equation*}
for almost every $v$, $\eta\in\RR^3$ and $I$, $I'>0$. Then $k_\kin$ appears as the standard kernel involved in the monatomic case, and, using \cite{Gr-63, BGPS-13}, we can state that $k_\kin$ satisfies a $L^2$-regularity property. 

Consequently, we only need to focus on $k_\intern$. In fact, a fine study allows to obtain pointwise estimates of $k_\intern$. More precisely, since $\zeta_2\in(-\delta,\delta)$, taking \eqref{e:hypbint} into account, we can write (see \cite[Lemma 4.5]{borsoni-phd}), for almost every $I$ and $I'\ge0$ such that $\max(I,I')\ge 1$, 
\begin{equation*}
    k_\intern(I,I') \lesssim e^{-\frac{|I'-I|}{4T_\intern}}~ [1+\max(I,I')]^{\delta-1-\zeta_2/2},
\end{equation*}
and, for almost every $I$ and $I'\ge0$ such that $\max(I,I')< 1$, 
\begin{equation*}
    k_\intern(I,I') \lesssim \int_{\max(I,I')}^2 U^{\zeta_2/2-1}\,\md U.
\end{equation*}
The previous bounds eventually lead to a $L^2$-regularity property for $k_\intern$. 

Then both regularity results on $k_\kin$ and $k_\intern$ straightforwardly imply the joint $L^2$-regularity for $k_3$, which finally allows to conclude on the compactness of $K_3$. 

\section{Ideas of compactness proofs for other polyatomic models}\label{s: main ideas other}

\subsection{Main ideas of the extension to a multi-component mixture}

Let us describe briefly how the ideas of Bernhoff \cite{Be-23b} can be extended to treat the case of a mixture \cite{Be-24a,Be-24b} (considering the difference between monatomic and polyatomic species).
The operator $K_{1}$ can still be treated in a quite straightforward way:
noting that (a multiple of) the exponent of the product $M_{i}M_{j\ast }$
 can now be recast as%
\begin{equation*}
m_{i}\frac{\left\vert v\right\vert ^{2}}{2}+m_{j}\frac{\left\vert v_{\ast
}\right\vert ^{2}}{2}+I \mathbf{1}_{i\in \mathcal{P}}+I_{\ast } \mathbf{1}_{j\in \mathcal{P}}=%
\frac{m_{i}+m_{j}}{2}\left\vert \frac{m_{i}v+m_{j}v_{\ast }}{m_{i}+m_{j}}%
\right\vert ^{2}+E_{ij}.
\end{equation*}%

For the operator $K_{3}$, as in the single species case, another parametrization is used, resulting in the following relations between velocities with $\xi\perp n$%
\begin{equation*}
\left\{ 
\begin{array}{c}
v_{\ast }=v^{\prime }+\xi+\chi _{-}n \\ 
v_{\ast }^{\prime }=v+\xi+\chi _{+}n%
\end{array}%
\right., \quad \text{where } \quad \chi _{\pm }=\frac{\Delta I_{ij}}{%
m_{i}\left\vert v-v^{\prime }\right\vert }\pm \frac{m_{i}-m_{j}}{2m_{j}}%
\left\vert v-v^{\prime }\right\vert .
\end{equation*}%
Unlike in the single species case (or species with the same 
mass), an additional term appears due to disparate masses $m_{i}\neq m_{j}$ (as in the monatomic case). Then (a multiple of) the exponent of the product $M_{j\ast
}M_{j\ast }^{\prime }$ can be recast as 
\begin{multline*}
\frac{\left\vert v_{\ast }\right\vert ^{2}}{2}+\frac{\left\vert v_{\ast
}^{\prime }\right\vert ^{2}}{2}+\frac{I_{\ast }+I_{\ast }^{\prime }}{m_{j}}\mathbf{1}_{j\in \mathcal{P}} \\
=\left\vert \frac{\left( v+v^{\prime }\right) _{\perp _{n}}}{2}%
+\xi\right\vert ^{2}+\frac{\left( m_{i}\left( \left\vert v^{\prime
}\right\vert ^{2}-\left\vert v\right\vert ^{2}\right) +2\Delta I_{ij}\right)
^{2}}{4m_{i}^{2}\left\vert v-v^{\prime }\right\vert ^{2}}+\frac{m_{i}^{2}}{%
m_j^{2}}\frac{\left\vert v-v^{\prime }\right\vert ^{2}}{4}+\frac{I_{\ast }+I_{\ast }^{\prime }}{m_{j}}\mathbf{1}_{j\in \mathcal{P}}\text{.}
\end{multline*}
Similar arguments as in the single species case can be applied for polyatomic species (or the collision operator for two monatomic
species). However, considering the collision operator when one (involved) species is monatomic and the other is polyatomic
causes some additional difficulties. The proof relies on using a suitable domain decomposition and careful estimates to proceed in the same way. 

Observe that for this operator $K_3$, the variables of the kernel are velocities and internal energies belonging to the same species. 
However, for the operator $K_2$, in the case of disparate masses, it is not the same anymore, and it is not obvious to find an analogous to the previous parametrization.
Nevertheless, the following inequality for (a multiple of) the exponent of the product $%
M_{i}^{\prime }M_{j\ast }$ is applied%
\begin{multline}\label{eq:ineq_mixt_continous}
m_{i}\frac{\left\vert v^{\prime }\right\vert ^{2}}{2}+m_{j}\frac{%
\left\vert v_{\ast }\right\vert ^{2}}{2}+I^{\prime }\mathbf{1}_{i\in \mathcal{P}}+I_{\ast }\mathbf{1}_{j\in \mathcal{P}} \\
\geq \dfrac{\left\vert m_{i}v-m_{j}v_{\ast }^{\prime }\right\vert
^{2}+m_{i}m_{j}\left\vert v-v_{\ast }^{\prime }\right\vert ^{2}+2\left(
m_{i}-m_{j}\right) \Delta I_{ij}}{2\left( \sqrt{m_{i}}+\sqrt{m_{j}}\right)
^{2}}+I^{\prime }\mathbf{1}_{i\in \mathcal{P}}+I_{\ast }\mathbf{1}_{j\in \mathcal{P}}\text{.}
\end{multline}
Then, with suitable domain decompositions, using different
parametrizations and convenient estimates, the kernel is shown to be $L^{2}$, and hence the operator $K_{2}$ to be a Hilbert-Schmidt integral operator.

Note that this inequality cannot be used in the case of equal masses $m_{i}=m_{j}$, but in this setting, the operator $K_{2}$ can be treated in a similar manner as the operator $K_{3}$ (shifting the roles of $w^{\prime }$ and $w_{\ast }^{\prime }$).

As for Brull, Shahine and Thieullen's approach, it can be extended in a straightforward way for the mixture case \cite{shahine-phd}. The same ideas are used as for the single species case, with different masses and different values of $\delta_i$, $\zeta_{ij}$ for $1\leq i,j \leq N$. As for the single species case, there is no additional difficulty to treat $K_3$, but due to the loss of symmetry, the treatment of the operators $K_2$ and $K_3$ lead to different integration assumptions in Hypothesis~\hypref{hyp:brull-mix}.

\subsection{Main ideas of the extension to discrete internal energy variables}
For single species, a similar strategy as for the case of a continuous
internal energy variable can be applied in the case of a discrete internal energy variable, without having to be concerned about any integral with respect to an internal energy. Instead, noting
that the set of internal energies is bounded above and below by%
\begin{equation*}
\min_{1\leq k\leq N_\intern} I^{(k)}>0\quad \text{ and } \quad 
\max_{1\leq k\leq N_{\intern}} I^{(k)}<+\infty \text{.}
\end{equation*}
This is also true for a multi-species mixture with discrete internal energy variables, with the following upper and lower bounds on the set of
internal energies %
\begin{equation*}
\min_{1\leq i\leq N}
\min_{1\leq k\leq N_{\intern,i}}I^{(k)}_i>0\quad \text{ and } \quad \max_{1\leq i\leq N}
\max_{1\leq k\leq N_{\intern,i}}I^{(k)}_i<+\infty \text{.}
\end{equation*}%
However, for the operator $K_{2}$ in the case of disparate masses, instead of applying the inequality \eqref{eq:ineq_mixt_continous} as in the continuous case, it is possible to use the following one for (a multiple of) the exponent of the product $M_{i}^{\prime }M_{j\ast }$ 
\begin{multline*}
m_{i}\frac{\left\vert v^{\prime }\right\vert ^{2}}{2}+m_{j}\frac{%
\left\vert v_{\ast }\right\vert ^{2}}{2}+I_{i}^{(k')}+I_{j}^{(\ell)}\\
\geq \left( \frac{\sqrt{m_{i}}-\sqrt{m_{j}}}{\sqrt{m_{i}}+\sqrt{m_{j}}}%
\right) ^{2}\frac{m_{i}\left\vert v\right\vert ^{2}+m_{j}\left\vert v_{\ast
}^{\prime }\right\vert ^{2}}{2}+\frac{\sqrt{m_{i}}-\sqrt{m_{j}}}{\sqrt{m_{i}}%
+\sqrt{m_{j}}}\Delta I_{ij}^{\left( k\ell,k^{\prime }\ell^{\prime }\right)
}+I_{i}^{(k')}+I_{j}^{(\ell)}\text{.}
\end{multline*}
Observe that this inequality reduces to \cite[Lemma 4.3]{BGPS-13} in the special case of a mixture of monatomic species (where all internal energies, and so all energy gaps as well, vanish) \cite{Be-23a}. With this form, the kernels for the components can be proved to be in $L^{2}((\RR^3)^2)$. However, unlike the monatomic case, different cases have to be considered, depending on the order of the relative speeds $%
\left\vert v-v_{\ast }\right\vert $ and $\left\vert v^{\prime }-v_{\ast
}^{\prime }\right\vert $ relatively to $1$ and $\left\vert v-v^{\prime
}\right\vert $. In the monatomic case though, the inequalities $\left\vert v-v_{\ast }\right\vert >\left\vert v-v^{\prime }\right\vert $ and $%
\left\vert v^{\prime }-v_{\ast }^{\prime }\right\vert >\left\vert
v-v^{\prime }\right\vert $ are always true. Different parametrizations are used
depending on the different cases.

\section{Discussion on the different hypotheses and results}\label{sec:discussion}

In this section, let us focus on the single species case, and on 
the following specific collision kernel
\begin{equation}\label{coll kernel fix}
B(v,v_*,I,I_*,r,R,\sigma) \propto E^{\zeta/2}.
\end{equation}
In order to compare the different hypotheses of Section \ref{s:compactness} to physically relevant cases, we shall first describe in the next subsection how physical values for $\delta$ and $\zeta$ can be obtained from experimental data.

\subsection{Range of physical values for the parameters}

In order to understand the range of physical values for $\delta$ and $\zeta$, let us determine the values of these parameters from measurements for some gases \cite{MPC-Dj-T-O, MPC-SS-non-poly}. 

The value of $\delta$ is related to experimental data on the specific heat $c_v$ through the equality \eqref{delta cv}, which holds under the polytropic assumption that $c_v$ remains constant. Practically, in \cite{MPC-Dj-T-O, MPC-SS-non-poly}, a gas is considered polytropic (calorically perfect) if its relative change of $c_v$ is below 5\% over a significant temperature range that starts with the room temperature of $300$~K. Examples of such gases are $\mbox{N}_2$,  $\mbox{O}_2$, CO, $\mbox{H}_2$.  

Next, the shear viscosity $\mu$ of a gas is known \cite{Ch-Cow} to depend on temperature $T$ through the power law
$$\frac{\mu(T)}{\mu(T_0)} = \left( \frac{T}{T_0} \right)^{{s}_{\textnormal{visc}}},$$
where $T_0$ is a reference temperature, and the viscosity exponent $s_{\textnormal{visc}}>0$ can be taken from \cite{Ch-Cow}, or extracted from experimental data through a fitting procedure. The evaluation of the Boltzmann collision operator \eqref{collisionoperator} with the collision kernel \eqref{coll kernel fix} under Chapman-Enskog asymptotics yields the same power-law behavior of the viscosity, with the exponent $\zeta/2-1$.
This analogy allows to exhibit values of $\zeta$. Hence, physical values of $\delta$ and $\zeta$ are presented in Table~\ref{Table:exp}, where those values are provided at different pressures, the ambient one ($1$\,bar) and a low one ($0.092$\,bar = $69$\,mmHg), and other values of $\zeta$ are also taken from \cite{Ch-Cow}.

\begin{table}[!htbp]
\caption{Physical values of $\delta$ and $\zeta$.}
	\begin{tabular}{| c || c | c | c | c | c |} \hline
	species &  	 \begin{tabular}{@{}c@{}}
						Polytropic temp. \\
						interval  (in K)\\
				\end{tabular}  
			& 
			\begin{tabular}{@{}c@{}}
					Pressure \\
						 (in bar)\\
				\end{tabular}  
			&  	 $\delta$
			& 	$\zeta$ 
			& $\zeta$ \cite{Ch-Cow}
\\ \hline \hline
\multirow{2}{*}{N$_2$} & \multirow{2}{*}{ $[300, 600]$} & 1  & 2.017 & 0.537 & \multirow{2}{*}{0.524}  \\   \cline{3-5}
&& 0.092 & 2.007 & 0.536 &  \\ \hline	 \hline	
\multirow{2}{*}{O$_2$} & \multirow{2}{*}{ $[300, 430]$}  & 1 & 2.080 & 0.443 & \multirow{2}{*}{0.454}  \\  \cline{3-5}
&& 0.092 & 2.070 & 0.441 & \\ \hline	\hline	
	\multirow{2}{*}{CO} & 	\multirow{2}{*}{$[300, 550]$} & 1 & 2.022 & 0.547  & \multirow{2}{*}{0.532}  \\ \cline{3-5}
& & 0.092 &	2.011 & 0.524  & \\ \hline \hline	
	\multirow{2}{*}{H$_2$} & 	\multirow{2}{*}{ $[300, 890]$} & 1 & 1.940  & 0.608  & \multirow{2}{*}{0.664}   \\  \cline{3-5}
&	& 0.092 &  1.939 & 0.608  &  \\ \hline
\end{tabular}
\label{Table:exp}
\end{table}	

\subsection{Comparison of the different hypotheses and link to physical values}

In the beginning of Section~\ref{s:compactness}, we list various hypotheses on the collision kernels in the form of upper and lower bounds, but formulated in, sometimes, very different ways, which are subsequently difficult to compare. Let us provide some comparisons in the framework of \eqref{coll kernel fix}.

On the first hand, considering $\delta \geq2$ and $0<\zeta \leq 2$ in \eqref{coll kernel fix} fits Hypothesis~\hypref{hyp:niclas}. In fact, the first term $E$ of the right-hand side of Hypothesis~\hypref{hyp:niclas} corresponds to the case $\zeta = 2$ in \eqref{coll kernel fix}. That same first term $E$ also bounds \eqref{coll kernel fix} for any $0<\zeta <2$ when $E\ge 1$. The case when $E\le 1$ is a bit more intricate to handle, but still fits \eqref{coll kernel fix} for $0<\zeta <2$ .  
\begin{remark}
By similar arguments, the restriction of Hypothesis \hypref{hyp:niclas-mix} to the single-species case holds for any $\delta \geq2$ and $-1<\zeta \leq 1$, which in combination with Hypothesis \hypref{hyp:niclas}  extends the applicable interval of $\zeta$ to $-1<\zeta \leq 2$. In fact, the upper limit of the interval can be extended, even if with doubtful physical meaning, to $-1<\zeta \leq \delta + 1$, see Remark~5 in \cite{Be-23b}. 
\end{remark}
On the other hand, the integrability condition from Hypothesis~\hypref{hyp:brull} couples values of $\delta$ and $\zeta$. Namely, in \eqref{coll kernel fix}, for $\zeta>-1$, the value of $\delta$ is restricted to $\delta >\max \left( 2,2+\zeta \right)$, due to the imposed integrability of $r^{\delta /2-2}(1-r)^{\delta -3-\zeta}$ over $(0,1)$. 

\medskip

It is thus interesting to compare the range of parameters considered in the hypotheses in Section~\ref{s:compactness} in view of the values from Table~\ref{Table:exp}. Let us first emphasize that Hypothesis~\hypref{hyp:brull} appears, at least for the choice \eqref{coll kernel fix} of the collision kernel, as unfitted for these gases, since the condition $\delta>\max(2,2+\zeta)$ is never satisfied. Second, Hypothesis~\hypref{hyp:niclas} is clearly satisfied by the nitrogen, oxygen and carbon monoxide values of $\delta$ and $\zeta$. But that is not the case for hydrogen. In a more general way, the case when $\delta<2$ is not covered yet by any existing results, while it does happen with some real gases. It is maybe possible to sharply adapt the assumptions on the parameters and the proofs so that the compactness results hold accordingly with respect to the physical values, but it is clearly beyond the scope of this review paper.

\section*{Acknowledgements}
The four authors acknowledge  the  support from COST Action CA18232 MAT-DYN-NET. 
The first author acknowledges travel grants by MAP5 at Université Paris Cité and SVeFUM, as well as the kind hospitality of MAP5.

\bibliographystyle{abbrv}
\bibliography{biblio}

\end{document}